\documentclass[12pt]{article}
\usepackage{amssymb}

\usepackage{mitpress}

\newtheorem{theorem}{Theorem}

\newtheorem{corollary}[theorem]{Corollary}

\newtheorem{definition}[theorem]{Definition}

\newtheorem{lemma}[theorem]{Lemma}

\newtheorem{proposition}[theorem]{Proposition}

\newdimen\dummy
\dummy=\oddsidemargin
\input{tcilatex}

\begin{document}

\title{Weak Omega Categories I}
\author{Carl A. Futia \\
5501 Keeney Street, Morton Grove, IL 60053\\
e-mail: topos8@aol.com}
\maketitle
\date{April 12, 2004}

\begin{abstract}
We develop a theory of weak omega categories that will be accessible to
anyone who is familiar with the language of categories and functors and who
has encountered the definition of a strict 2-category.

The most remarkable feature of this theory is its simplicity. We build upon
an idea due to Jacques Penon by defining a weak omega category to be a span
of omega magmas with certain properties. (An omega magma is a reflexive,
globular set with a system of partially defined, binary composition
operations which respects the globular structure.)

Categories, bicategories, strict omega categories and Penon's weak omega
categories are all instances of our weak omega categories. We offer a
heuristic argument to justify the claim that Batanin's weak omega categories
also fit into our framework.

We show that the Baez-Dolan stabilization hypothesis is a direct consequence
of our definition of weak omega categories.

We define a natural notion of a pseudo-functor between weak omega categories
and show that it includes the classical notion of a homomorphism between
bicategories. In any weak omega category the operation of composition with a
fixed 1-cell defines such a pseudo-functor.

Finally, we define a notion of weak equivalence between weak omega
categories which generalizes the standard definition of an equivalence
between ordinary categories.
\end{abstract}

\subsection{\protect\bigskip Introduction}

This paper begins the development of a theory of weak, higher dimensional
categories which parallels closely the familiar theory of ordinary (1
dimensional ) categories [15]. It should be accessible to anyone comfortable
with the language of categories and functors. In particular, concepts such
as operads and monads play no role in the basic definitions. In the sequel,
Part II, we shall build upon the foundation laid down in Part I to develop
some of the more technical aspects of the theory: a notion of weighted
limits and the construction of the weak omega category of small weak omega
categories. In Part II we shall also construct a functor from any Quillen
model category to this omega category of small omega categories and
construct a weak omega category whose n-dimensional morphisms are
n-dimensional cobordisms with corners.

A category can be defined as a directed graph with a partially defined,
binary composition law that satisfies additional axioms: associativity of
composition and the existence of right and left identity for each vertex (or
object). Our weak omega categories are defined in a similar spirit.

We start with an omega graph (a reflexive, globular set) together with a
system of partially defined, binary composition laws that respects the graph
structure. Such an omega graph together with its system of composition laws
is called an omega magma.

An omega magma is a strict omega category if its composition laws satisfy
the higher dimensional generalizations of the associative, identity and
interchange laws exactly, i.e. these axioms hold as equations between
elements of the omega magma. A strict omega category all of whose cells
above dimension 1 are identities is just an ordinary category. The standard
example of a strict 2-category [2,7](all cells above dimension 2 are
identities) is the 2-category whose objects are small categories, whose
morphisms are functors and whose 2-dimensional arrows are natural
transformations between functors.

An omega magma is a weak omega category if its composition laws satisfy the
higher dimensional generalizations of the associative, identity and
interchange laws in a ``relaxed'' way: the laws are required to hold only
``up to an equivalence''. An equivalence in a weak omega category is a
higher dimensional arrow in the category that behaves like a homotopy
equivalence in homotopy theory. It is generally something less than an
isomorphism but still preserves enough structure to behave like one for
categorical purposes. The simplest example is already familiar to the
reader: any functor inducing a categorical equivalence between a pair of
ordinary categories is a 1-dimensional equivalence in this more general
sense (at least if we admit the axiom of choice); such a functor need not be
an isomorphism of categories.

Any bicategory [6,7,12,14,15,16] is a weak 2-category. Monoidal categories
[8,15 ] are bicategories having only one object and are the most familiar
instances of weak 2-categories.

The crux of any theory of weak omega categories is its method for expressing
mathematically the coherence conditions that assert that the desired
categorical laws hold ``up to equivalence''. The device we have chosen for
this purpose first appeared implicitly in the work of Jacques Penon [18]. It
traces its roots to the standard coherence theorem for bicategories [12,14].
The latter result asserts that any bicategory can be embedded in a strict
2-category via a functor which preserves the 1-dimensional composition law
only up to isomorphism. Penon's wonderful idea was simply to turn the
conclusion of this coherence theorem into a definition.

Of course some subtlety must be involved here. Category theorists have known
for some time [11] that there are weak 3-categories which cannot be embedded
into a strict 3-category. Thus it is not possible to define a weak omega
category as an omega magma that admits an appropriate, structure preserving
embedding into a strict omega category.

Penon circumvented this problem by redefining a morphism from one omega
magma to another to be a span of omega magma homomorphisms, i.e. a diagram
of the form $X\leftarrow Z\rightarrow Y$ in the category of omega magmas. We
call this a span from $X$ to $Y$ and call $X$ the domain and $Y$ the
codomain of the span. One can then define a weak omega category to be an
omega magma that is the domain of a span to a strict omega category. To
avoid a trivial theory one must of course impose some additional conditions
on this span.

We should warn the reader that he will not find this informal explanation
anywhere in Penon's paper [18]. Nonetheless it lies just beneath the surface
of his work and soon becomes apparent once one attempts to unravel his
definition of weak omega categories (which he called ''prolixes'').

The approach to weak omega categories taken here differs in several other
ways from Penon's.

We highlight the idea that a weak omega category is first of all an omega
magma: an omega graph with a system of partially defined, binary composition
operations. The coherence conditions that make the omega magma a weak omega
category are expressed by a particular span from the omega magma to a strict
omega category. This span must satisfy certain simple axioms. In Penon's
theory omega magmas are incidental, simply stepping stones on the way to the
construction of a monad on the category of reflexive, globular sets whose
algebras are his weak omega categories.

We might add that our emphasis on omega magmas and their systems of binary
compositions also distinguishes this work from theories of higher
dimensional categories that define them as algebras for higher dimensional
operads [3,5,9,13,14,17].

A much more important difference is our willingness to entirely remove
certain restrictions Penon places upon the spans which define his weak omega
categories . If $X\leftarrow Z\rightarrow Y$ is a span defining $X$ as a
weak omega category, Penon requires that the strict omega category $Y$ be
freely generated by the omega graph underlying $X$. In addition the omega
magma $Z$ is required to be freely constructed (via an adjunction) from the
free omega magma generated by the same omega graph, viz. the one underlying $%
X$. (The omega magma $Z$ Penon calls a ''stretching'' of the strict omega
category $Y$.) These last restrictions are inherent in Penon's construction
of the monad whose algebras are his weak omega categories.

It soon became apparent to us that these requirements render Penon's theory
inflexible, making constructions and proofs hard and obscure in situations
where they should be easy and transparent. We therefore have chosen to allow
the strict omega category $Y$ and the omega magma $Z$ to be completely
arbitrary while retaining certain axioms on the span from $X$ to $Y$. The
price we pay for this is that our weak omega categories are no longer the
algebras for some suitable monad on omega graphs. Still, the benefits of
this generality far outweigh its costs.

The biggest benefit lies in the great simplicity of the resulting theory.
This allows us to explore territory as yet inaccessible to other theories of
weak omega categories [3,5,9,13,14,17,18,19] (see the next subsection for
details).

Many mathematicians have offered encouragement, insight and patient answers
to our often benighted questions during the course of this research. For
their generous assistance we would like to thank John Baez, Michael Batanin,
Clemens Berger, Ronnie Brown, Eugenia Cheng, John Duskin, Paul Goerss, Peter
Johnstone, G. Max Kelly, Steven Lack, Tom Leinster, Saunders MacLane, Peter
May, Tim Porter, Charles Rezk, Steven Schnauel, Ross Street, Earl Taft, Mark
Weber, Noson Yanofsky and David Yetter.

We would especially like to thank Michael Batanin for his kind encouragement
during the course of this research and for his willingness to answer
technical questions about his operadic methods when these arose. Moreover,
his insight that composites in strict omega categories can be effectively
represented and manipulated using his language of level trees was an
absolutely essential foundation for our own effort to understand weak omega
categories.

\subsection{Contents of this paper}

Section 1 begins with the definitions of globular sets, omega-graphs, omega
magmas and strict omega categories. We discuss briefly the concept of a
locally presentable category and the notion of an essentially algebraic
theory. The category of models of any essentially algebraic theory is a
locally presentable category. Locally presentable categories have many
convenient properties which facilitate the construction of adjunctions.

Section 2 presents the definition of a weak omega category$\ \mathbf{X}$ as
a span of omega magmas that satisfies certain simple axioms. The domain of
such a span is the underlying omega magma of the weak omega category. We
introduce the notion of a Penon map between omega magmas. Such maps appear
as the codomain ``leg'' of the span defining a weak omega category and they
structure the coherence data for the span. We define an obvious notion of an
omega functor between weak omega categories and show that the resulting
category, \textit{Omega\_Cat, }of small, weak omega categories and omega
functors is the category of models of an essentially algebraic theory.

In Section 3 we begin to justify our definition of weak omega categories by
showing that \textit{Cat}, the category of small categories and functors, is
a retract of the full subcategory of \textit{Omega\_Cat} containing those
objects which are \textit{1-skeletal}, i.e. those objects whose underlying
omega magma (the domain of the defining span) has only identity cells above

\noindent dimension 1. We also show that \textit{Strict\_Cat}, the category
of strict omega categories and strict omega functors is in fact isomorphic
to a full subcategory of \textit{Omega\_Cat.}

We continue justifying our definition in Section 4 by considering
bicategories. We show that \textit{Bicat}, the category whose objects are
bicategories and whose morphisms are strong homomorphisms between
bicategories (preserving operations and identities ``on the nose'') is a
retract of the full subcategory of \textit{Omega\_Cat} whose objects are 
\textit{2-skeletal}.

Section 5 examines the relationship between two other definitions of weak
omega categories and our own. We show that the category of Penon's weak
omega categories and omega functors, \textit{Prolixe}, is a retract of a
certain full subcategory of our \textit{Omega\_Cat. }We also offer an
informal argument to support our contention that \textit{Batanin\_Cat}, the
category whose objects are Batanin's weak omega categories (i.e. the
algebras for the initial contractible, higher dimensional operad with a
system of compositions) are instances of our weak omega categories.

In Section 6 we point out one immediate consequence of our definition. For
each $k\geq 1$ one can define a full subcategory \textit{nCat}$_{k}$ of 
\textit{Omega\_Cat }in which an object is a weak omega category whose three
defining omega magmas each has only a single cell in every dimension $\leq
k-1$. \ Objects of \textit{nCat}$_{k}$ are called k-tuply monoidal weak $%
(n+k)$-categories by Baez and Dolan [4\.{]}. \ We define an obvious pair of
functors in opposite directions connecting \textit{nCat}$_{k}$ and \textit{%
nCat}$_{k+1}$ whenever $k\geq 2$ and observe that each composite functor is
the identity functor. These functors simply shift all the data defining weak
omega categories and omega functors by a single dimension. In this sense the
Baez-Dolan stabilization conjecture [4] is a simple consequence of our
definition of weak omega category.

Section 7 discusses several methods for constructing weak omega categories
from given ones and for recognizing that a given omega magma can be given
the structure of a weak omega category.

\textit{Omega\_Cat} is the category of models of an essentially algebraic
theory. As such it is locally (finitely) presentable and thus complete and
cocomplete with respect to ordinary conical limits and colimits. We show
that the left adjoint half of Penon's monad [18] can be used to define a
functor from \textit{Omega\_Graph} to \textit{Omega\_Cat }and thus gives a
way of functorially associating with any omega graph a weak omega category.
Moreover, Penon's methods also show that to any morphism from an omega graph
to a strict omega category one can associate functorially a weak omega
category whose defining span has as its codomain the given strict omega
category.

If $\mathbf{X}$ is a weak omega category and $a,b$ a pair of cells of
dimension $i-1\geq 0$ we construct the weak omega category $\mathbf{X}(a,b)$
whose objects are $i$ cells of $\mathbf{X}$ with domain $a$ and codomain $b$.

Section 7 ends with a definition of \textit{categorical equivalence relations%
} and offers a condition sufficient to guarantee that an omega magma
equipped with such a relation is the underlying magma of a weak omega
category.

In Section 8 we define a natural notion of omega pseudo-functor between weak
omega categories. Every omega functor is an omega pseudo-functor. We show
that omega pseudo-functors compose in the obvious way and that therefore
there is an ordinary category \textit{PF\_Omega\_Cat} whose objects are weak
omega categories and whose morphisms are omega pseudo-functors. If $\mathbf{X%
}$ is a weak omega category we show that the operation of composition with a
fixed object of $\mathbf{X}(b,c)$ is a pseudo-functor from $\mathbf{X}(a,b)$
to $\mathbf{X}(a,c)$ for any triple of $i-1$ cells $\ a,b,c$.

We define the notion of a \textit{proper} homorphism between classical
bicategories and show that any proper homomorphism between bicategories is
an omega pseudo-functor. We know of no examples of classical homomorphisms
that are not proper but in principle such examples could exist. We also show
that any omega pseudo-functor between 2-skeletal, weak omega categories is a
proper homomorphism between classical bicategories.

We conclude in Section 9 by defining the notions of \textit{weak equivalence 
}and of \textit{omega equivalence }between weak omega categories. These
definitions are the foundation of a theory of weighted limits which shall be
developed in part II of this work. The basic idea is to define a suitable
\textquotedblleft components\textquotedblright\ functor, $\Pi $, from a
certain category \textit{Tame\_Omega\_Cat} to \textit{Set}, the category of
small sets. The category \textit{Tame\_Omega\_Cat} has as its objects small,
weak omega categories and as its morphisms those omega pseudo functors which
we shall call \textit{tame.} An omega pseudo-functor is tame when it
preserves a class of arrows we call internal equivalences, much as a functor
preserves isomorphisms and as a homomorphism between bicategories preserves
1-cell equivalences. The class of tame omega pseudo-functors includes all
omega functors between weak omega categories and all omega pseudo-functors
between weak n-categories for n finite.

We say that an omega pseudo-functor $\mathbf{F:X\rightarrow Y}$ is a weak
equivalence if $\Pi (\mathbf{F})$ is an isomorphism and if $\Pi (\mathbf{F}%
(a,b))$ is also an isomorphism for all pairs of $i-1$ cells $a,b$ where 
\[
\mathbf{F}(a,b):\mathbf{X}(a,b)\rightarrow \mathbf{Y}(Fa,Fb) 
\]
is the restriction of $\mathbf{F}$ to\textbf{\ }$\mathbf{X}(a,b)$. It is
easily seen that two weak omega categories which are 1-skeletal (and hence
are ordinary categories) are weakly equivalent in this sense if and only if
they are equivalent as ordinary categories and the appropriately modified
assertion also holds for 2-skeletal, weak omega categories..

\section{Recollections}

\subsection{\protect\bigskip Omega magmas and strict omega categories}

We begin by recalling some definitions which are essential to our theory of
higher dimensional categories. \ These are all standard and can be found for
example in [5,14] (with some differences in notation). The only one which
may be unfamiliar is the notion of an omega magma (compare [14,18]). An
omega magma can be understood as an object that would be a strict omega
category but for the failure of its composition laws to satisfy the
associative, interchange and identity laws which must hold in any strict
omega category.\vspace{0in}

\begin{definition}
A \textbf{globular set }$X$ is a sequence of sets $(X_{i},i\geq 0)$ together
with functions $dom_{i-1}^{i},cod_{i-1}^{i}:X_{i}\rightarrow X_{i-1}$
defined for $i>0$ t$\bigskip $hat satisfy the so-called \textbf{globular
relations}: 
\[
dom_{i-1}^{i}\circ dom_{i}^{i+1}=dom_{i-1}^{i}\circ cod_{i}^{i+1} 
\]
\[
cod_{i-1}^{i}\circ dom_{i}^{i+1}=cod_{i-1}^{i}\circ cod_{i}^{i+1} 
\]
for all integers $i>0$
\end{definition}

The elements of the set $X_{i}$ are called the \textit{cells} of $X$ of
dimension $i$. The function $dom_{i-1}^{i}$ assigns to each $i$-cell its
domain and the function $cod_{i-1}^{i}$ assigns to each $i$-cell its
codomain. The terminology is intended to evoke the idea that an $i$-cell of $%
X$ is a kind of \textquotedblleft $i$-dimensional
morphism\textquotedblright\ whose domain and codomain are both
\textquotedblleft $(i-1)$-dimensional morphisms\textquotedblright .

For $j-i>0$ we denote the $(j-i)$-fold composite 
\[
dom_{i}^{i+1}\circ dom_{i+1}^{i+2}\circ ...\circ dom_{j-1}^{j} 
\]
by $dom_{i}^{j}$.

There is a category, \textit{Glob\_Set,} in which an object is a small
globular set and in which a morphism $F:X\rightarrow Y$ is a sequence of
functions $(F_{i}:X_{i}\rightarrow Y_{i})$ which commute with the domain and
codomain functions of $X$ and $Y$.

Any directed graph $G$ determines a globular set by defining $G_{0}$ to be
the set of vertices and $G_{i}$ to be the set of directed edges for $i\geq 1$%
. For

\noindent $i>1$ $dom_{i-1}^{i}$ is the identity map and for $i=1$ this
function assigns to a directed edge the source vertex for that edge. \ For $%
i>1$ $\ cod_{i-1}^{i}$ is also the identity map and for $i=1$ this function
assigns to a directed edge its target vertex. Thus to any category there is
a globular set \ determined by the category's underlying directed graph
whose vertices are the objects of the category and whose directed edges are
the morphisms.

The directed graph associated with a category has additional structure, the
directed edges corresponding to the identity morphisms, that isn't part of
the definition of a globular set. Globular sets with this additional
structure we shall call \textit{omega graphs}.

\begin{definition}
An \textbf{omega graph} $X$ is a globular set together with a sequence of 
\textbf{identity functions} $id_{i}^{i-1}:X_{i-1}\rightarrow X_{i}$ \ for $%
i>0$ which satisfy the following relations: 
\[
dom_{i-1}^{i}\circ id_{i}^{i-1}(x)=x 
\]
\[
cod_{i-1}^{i}\circ id_{i}^{i-1}(x)=x 
\]
for all $x\in X_{i-1}$. An $i$-cell of $X$ in the image of $id_{i}^{i-1}$ is
called an \textbf{identity cell} of $X$. The $(j-i)$-fold iterated composite 
$id_{j}^{i}$ of identity functions is defined by the equation 
\[
id_{j}^{i}=id_{i+1}^{i}\circ id_{i+2}^{i+1}\circ ...\circ id_{j}^{j-1} 
\]
\end{definition}

In the literature one often encounters omega graphs under a different
name:reflexive, globular sets. We prefer our own terminology both because of
its economy and because its use avoids the risk of confusing globular sets
with reflexive, globular sets.

There is a category, \textit{Omega\_Graph}, in which an object is a small
omega graph and a morphism is a morphism of the underlying globular sets
which also commutes with the identity cell functions.

We shall often have occasion to consider omega graphs in which all cells
above a certain dimension are identity cells.

\begin{definition}
An omega graph $X$ is\textbf{\ }$n$\textbf{-skeletal} if every cell of
dimension $>n$ is an identity cell.
\end{definition}

Our theory of weak omega categories is built around the notion of an omega
graph with a compatible family of partially defined, binary composition laws.

\begin{definition}
An \textbf{omega magma }is an omega graph together with a sequence of
ternary relations 
\[
Comp_{i}^{j}\subseteq X_{j}\times X_{j}\times X_{j}\text{ \ ,}0\leq i<j 
\]
These relations are called \textbf{composition relations} and must satisfy
the following axioms:

\begin{enumerate}
\item  if $(a,b,c)\in Comp_{i}^{j}$ then $cod_{i}^{j}(a)=dom_{i}^{j}(b)$

\item  if $(a,b,c)\in Comp_{i}^{j}$ and $(a,b,d)\in Comp_{i}^{j}$ then $c=d$

\item  if $a,b\in X_{j}$ and $cod_{i}^{j}(a)=dom_{i}^{j}(b)$ then $\exists c 
$ and $(a,b,c)\in Comp_{i}^{j}$

\item  if $i+1=j$ and $(a,b,c)\in Comp_{i}^{j}$ then $dom_{i}^{j}(c)=$ $%
dom_{i}^{j}(a)$ and $cod_{i}^{j}(c)=cod_{i}^{j}(b)$

\item  if $i+1<j$ and $(a,b,c)\in Comp_{i}^{j}$ then 
\[
(dom_{j-1}^{j}(a),dom_{j-1}^{j}(b),dom_{j-1}^{j}(c))\in Comp_{i}^{j-1} 
\]
and 
\[
(cod_{j-1}^{j}(a),cod_{j-1}^{j}(b),cod_{j-1}^{j}(c))\in Comp_{i}^{j-1} 
\]
\end{enumerate}
\end{definition}

Each of the ternary relations which is part of the structure of an omega
magma $X$ is a partially defined, binary composition law. It is usually more
convenient to use the ``infix'' notation to denote composite cells in an
omega magma. Thus if $a,b\in X_{j}$ and $cod_{i}^{j}(a)=dom_{i}^{j}(b)$ we
shall write $a\bigcirc _{i}^{j}b$ to denote the unique $j$-cell $c$ for
which $(a,b,c)\in Comp_{i}^{j}$. \ Moreover, any appearance of the notation $%
a\bigcirc _{i}^{j}b$ will be taken to imply that $a,b\in X_{j}$ and $%
cod_{i}^{j}(a)=dom_{i}^{j}(b)$. This convention will save much tedious
repetition of obvious hypotheses. The partial operation $\bigcirc _{i}^{j}$
will sometimes be referred to as composition of $j$-cells over $i$-cells.

A category $C$ can be thought of as a 1-skeletal omega magma whose only non
trivial composition law is the relation $Comp_{0}^{1}$. If $(f,g,h)\in
Comp_{0}^{1}$ then the codomain of the morphism $f$ coincides with the
domain of the morphism $g$ and the composite $g\circ f$ is equal to $h$.

Here we point out a potential source of confusion. \ In any category it is
conventional to write the composite morphism of the diagram 
\[
\begin{array}{ccccc}
& f &  & g &  \\ 
X & \rightarrow & Y & \rightarrow & Z
\end{array}
\]
as $g\circ f$. On the other hand, if we are thinking of the category as a
1-skeletal omega magma the same element will be written as $f\bigcirc
_{i}^{j}g$. Note the reversal of order!

There is a category, \textit{Omega\_Magma}, in which an object is a small
omega magma and a morphism is a morphism of underlying omega graphs which
commutes with the partially defined composition operations.

Omega magmas can be quite complicated objects since their composition laws
satisfy no axioms aside from those which ensure compatibility with the omega
graph structure. However a consideration of low dimensional examples helps
to give some insight into the significance of their defining axioms and
shows how omega magmas differ from categories.

As a simple first example let's consider a 1-skeletal omega magma $X$ which
has only a single $0$-cell denoted by $\ast $ and which is generated in
dimension $1$ by a single $1$-cell $f$ together with the identity cell $%
id_{1}^{0}(\ast )$. Then the set $X_{1}$ of $1$-cells is just the set of all
possible ways of inserting parentheses into a finite sequence consisting of
repetitions of the symbols $f$ and $id_{1}^{0}(\ast )$ so as to represent a
meaningful sequence of binary compositions yielding a single cell. The
operation $\bigcirc _{0}^{1}$ is just concatenation of parenthesized
sequences. It is important to remember that while we call $id_{1}^{0}(\ast )$
an identity cell the two elements of $X_{1}$ denoted by $id_{1}^{0}(\ast )$ $%
\bigcirc _{0}^{1}f$ and $f$ (for example) are \textbf{distinct}.

Here is a slightly more general example which illustrates a connection
between omega magmas and labelled, rooted binary trees. Let $G$ be a
directed graph and let $X$ denote the 1-skeletal omega magma freely
generated by the graph $G$. The set $X_{0}$ is just the set of vertices of $%
G $. The set $X_{1}$ is a set of ordered pairs. The first coordinate of such
a pair is a finite sequence of composable edges of $G$ which will generally
include some identity edges corresponding to elements of $X_{0}$. The second
component is a specific choice of a way of inserting parentheses into the
sequence that is the first coordinate so as to represent a sequence of
binary compositions yielding a single cell. The operation $\bigcirc _{0}^{1}$
is again the obvious operation derived from concatenation of compatible
sequences of edges.

There is another way to represent the elements of $X_{1}$ in this example.
Every such element corresponds to a unique rooted, binary tree whose leaves
are linearly ordered and labelled. The labels of the leaves are just the
edges appearing in the sequence that is the first coordinate of the chosen
element of $X_{1}$. Distinct elements of $X_{1}$are represented by distinct
labelled trees. However not every labelled tree corresponds to an element in 
$X_{1}$; for this to be so the labels must represent a sequence of elements
that are compatible for composition in the order dictated by the tree
structure. In this representation the operation $\bigcirc _{0}^{1}$ places
side by side the trees representing the elements being composed and then
joins them by adding a new root.

We can now define the notion of a strict omega category. These will be seen
to be omega magmas whose partial operations are well behaved in that they
satisfy the associative, interchange and identity laws.

\begin{definition}
A \textbf{strict omega category} is an omega magma $X$ satisfying the
following axioms:

\begin{enumerate}
\item  \textbf{Associativity}. $(a\bigcirc _{i}^{j}b)\bigcirc
_{i}^{j}c=a\bigcirc _{i}^{j}(b\bigcirc _{i}^{j}c)$

\item  \textbf{Interchange.} If $i<j<k$ then 
\[
(a\bigcirc _{j}^{k}b)\bigcirc _{i}^{k}(c\bigcirc _{j}^{k}d)=(a\bigcirc
_{i}^{k}c)\bigcirc _{j}^{k}(b\bigcirc _{i}^{k}d) 
\]

\item  \textbf{Identity}. \ $a=id_{j}^{i}(dom_{i}^{j}a)\bigcirc
_{i}^{j}a=a\bigcirc _{i}^{j}id_{j}^{i}(cod_{i}^{j}a)$

\item  I\textbf{identity Interchange.} \ If\textbf{\ } $i<j$ then 
\[
id_{j+1}^{j}(a)\bigcirc _{i}^{j+1}id_{j+1}^{j}(b)=id_{j+1}^{j}(a\bigcirc
_{i}^{j}b) 
\]
\end{enumerate}
\end{definition}

There is a category, \textit{Strict\_Category}, in which an object is a
small, strict omega category and a morphism is a morphism of underlying
omega magmas. This is clearly a full subcategory of \textit{Omega\_Magma}.

Any strict 2-category in the usual sense is in an obvious way a 2-skeletal,
strict omega category. The standard example of such a strict 2-category is
the 2-category in which objects are small categories, morphisms are functors
and 2-cells natural transformations between functors.

Given a strict omega category $X$ and a triple of indices $i<j<k$ one sees
immediately that the sets $X_{i},X_{j}$ and $X_{k}$ are the zero, one and
two-cells respectively of a strict 2-category. In particular, for any pair
of indices $i<j$ the sets $X_{i}$ and $X_{j}$ are the objects and morphisms
respectively of an ordinary category in which the composition law is $%
\bigcirc _{i}^{j}$.

\subsection{Locally finitely presentable categories and essentially
algebraic theories}

In the Introduction to this paper we asserted that the reader needed only a
familiarity with the basic language of categories and functors to understand
our definition of weak omega category. Clearly the topics of this subsection
require more expertise.\ However, the reader who feels that locally
presentable categories and essentially algebraic theories are more technical
baggage than he or she wishes to carry may safely skip this subsection and
move on to Section 2. The material we shall now discuss is used only to
construct certain adjunctions between categories. Nothing is lost if the
reader is willing to take the existence of these adjunctions on faith.

Our basic references for locally presentable categories are [1] and [8]. For
a treatment of essentially algebraic theories the reader may consult [1].

\begin{definition}
An object $a$ of a category $A$ is \textbf{finitely presentable} if the
functor $A(a,-)$ preserves directed colimits.
\end{definition}

The reader will recall that a directed colimit is a colimit over a diagram
that is a directed poset. A poset is directed if each pair of elements has
an upper bound.

Perhaps the most familiar example of a finitely presentable object in a
category is a finitely presentable group (finite number of generators and
relations) in the category of groups.

\begin{definition}
A category $A$ is \textbf{locally finitely presentable} if it is cocomplete
and if there is a set $K$ of finitely presentable objects such that each
object $a$ of $A$ is a directed colimit of objects in $K$.
\end{definition}

Thus the category of groups is locally finitely presentable with the set $K$
consisting of one representative from each isomorphism class of finitely
presented groups.

There is a simple criterion for a functor between locally finitely
presentable categories to be a right adjoint and this is the reason for our
interest in such categories.

\begin{theorem}
([1] theorem 1.66) A functor between locally finitely presentable categories
is a right adjoint if and only if it preserves limits and directed colimits.
\end{theorem}

Occasionally we shall have reason to consider subcategories of locally
finitely presentable categories and shall wish to prove that they are also
locally finitely presentable. The following result will be useful in this
regard.

\begin{theorem}
([1] theorem 1.20) A category is locally finitely presentable if and only if
it is cocomplete and has a strong generator consisting of finitely
presentable objects.
\end{theorem}

We recall that a \textit{strong generator} for a category $A$ is a set 
\textsf{G} of objects with the following property. If $b$ and $c$ are
objects of $A$ and $b$ is a proper subobject $c$ then there is an object $%
a\in \mathsf{G}$ and a morphism $f:a\rightarrow c$ that does not factor
through $b$.

The preceding result asserts that every locally finitely presentable
category is cocomplete, \textit{but in fact such a category is also complete}
([1] theorem 1.28).

One way to show that a category is locally finitely presentable is to
exhibit it as the category of models of an essentially algebraic theory in
which all operations have finite arities, i.e. each depends only on a finite
number of arguments. All of the categories defined earlier in this section
can be seen to be locally finitely presentable for this reason. We next
informally explain the concept of an essentially algebraic theory (see [1]
chapter 3.D for the full story).

One defines an \textit{essentially algebraic theory} by starting with a set
of \textit{sorts}. For example, the theory of globular sets is essentially
algebraic and its $i^{th}$ sort is the sort of $i$-cells. Thus this theory
has one sort for each non-negative integer. The theory of omega magmas is
also essentially algebraic and it has the same set of sorts as the theory of
globular sets.

The second ingredient of an essentially algebraic theory is a set of \textit{%
total operations}. The total operations of the theory of globular sets are
the operations $dom_{i}^{j}$ and $\ cod_{i}^{j}$ for all pairs of integers $%
i<j$. The total operations of the theory of omega magmas are these plus the
operations $id_{j}^{i}$ for all pairs of integers $i<j$.

Essentially algebraic theories are distinguished from algebraic theories
because they allow some operations to be only partially defined. All
operations in the theory of globular sets are total and so this essentially
algebraic theory is actually an algebraic theory. On the other hand the
theory of omega magmas has \textit{partial operations} $\bigcirc _{i}^{j}$%
defined only for those ordered pairs of $j$-cells $(a,b)$ for which $%
cod_{i}^{j}a=dom_{i}^{j}b$. The important point here is that the domain of
definition for each partial operation is defined by equations involving only
the total operations.

The final ingredient needed to define an essentially algebraic theory is a
list of equations between operations (total and partial) which records the 
\textit{axioms} of the theory.

A \textit{model} of an essentially algebraic theory is determined by
assigning to each of the theory's sorts a small set and to each total and
partial operation an appropriate function between the (subsets of products
of) sets assigned to sorts. These functions must satisfy the axioms of the
theory. A \textit{morphism between models }of the theory is simply a
collections of functions, one for each sort, that maps the set corresponding
to a sort in one model to the set corresponding to the same sort in the
other model. \ These maps must of course commute with all the total and
partial operations.

There is thus defined a category in which an object is a model of the
essentially algebraic theory and a morphism is a morphism between models as
described above. This category is called the \textit{category of models of
the theory}.

The important property of the category of models of an essentially algebraic
theory is described by the following result.

\begin{theorem}
([1] theorem 3.36) A category is locally finitely presentable if and only if
it is equivalent to the category of models of an essentially algebraic
theory in which all arities are finite, each partial operation has domain
defined by a finite number of equations, each such equation involves a
finite number of variables and each axiom of the theory involves only a
finite number of variables.
\end{theorem}

It should be obvious from the definitions that \textit{Glob\_Set,
Omega\_Graph, Omega\_Magma} and \textit{Strict\_Category} are categories of
models of essentially algebraic theories and are moreover locally finitely
presentable by the preceding result.

\subsection{Some useful adjunctions}

The definitions in section 1.1 show that there is a sequence of categories
and forgetful functors 
\[
\begin{array}{ccccccc}
& U_{S}^{G} &  & U_{G}^{M} &  & U_{M}^{C} &  \\ 
\mathit{Glob\_Set} & \leftarrow & \mathit{Omega\_Graph} & \leftarrow & 
\mathit{Omega\_Magma} & \leftarrow & \mathit{Strict\_Category}
\end{array}
\]

Each of the categories in this diagram is the category of models of an
essentially algebraic theory and by theorem 10 is a locally finitely
presentable category. In particular each of these categories is complete and
cocomplete. The forgetful functors are easy to describe. $U_{M}^{C}$ forgets
the axioms of a strict omega category,\ $U_{G}^{M}$ forgets the composition
laws of an omega magma and $U_{S}^{G}$ forgets the identity functions of an
omega graph. It is equally easy to see that each of these functors creates
both limits and directed colimits (see [15, p.109] for the definition of
creating limits) and therefore preserves such limits and colimits.
Consequently theorem 8 assures us that each of these forgetful functors has
a left adjoint. Thus we obtain a diagram 
\[
\begin{array}{ccccccc}
& U_{S}^{G} &  & U_{G}^{M} &  & U_{M}^{C} &  \\ 
\mathit{Glob\_Set} & \leftrightarrows & \mathit{Omega\_Graph} & 
\leftrightarrows & \mathit{Omega\_Magma} & \leftrightarrows & \mathit{%
Strict\_Category} \\ 
& L_{G}^{S} &  & L_{M}^{G} &  & L_{C}^{M} & 
\end{array}
\]

These left adjoints have explicit descriptions. The composite functor $%
L_{C}^{S}=L_{C}^{M}\circ L_{M}^{G}\circ L_{G}^{S}$ \ is isomorphic to the
free omega category functor constructed by \ Batanin [5]. The functor $%
L_{G}^{S}$ freely adjoins to a globular set the required identity elements.
The composite $L_{C}^{G}=L_{C}^{M}\circ L_{M}^{G}$ is isomorphic to the free
omega category functor (generated by an omega graph) constructed by Penon
[18]. The functor $L_{M}^{G}$ is implicit in Penon's construction of his
``stretching'' omega magmas [18]. Finally, the functor $L_{M}^{G}$ assigns
to an omega magma the coequalizer obtained by imposing the relations which
must hold in any strict omega category.

These left adjoints will prove useful in analyzing the properties of our
definition of weak omega categories and for constructing examples.

\section{Weak omega categories}

In this section we offer our definition of weak omega categories. In
Sections 3,4 and 5 we shall explain how various types of higher dimensional
categories that have appeared in the literature are all instances of the
kind of weak omega category defined here.

Our first task is to define the concept of a \textit{bridge relation}. Such
relations can be thought of as carriers for the coherence data that defines
an omega category structure on an omega magma.

First we establish some notation. Let $X$ be an globular set and $a,b\in
X_{i}$. We say that $a$ and $b$ are \textit{parallel }and write $a\parallel
b $ if $dom_{i-1}^{i}a=dom_{i-1}^{i}b$ and $cod_{i-1}^{i}a=cod_{i-1}^{i}b$.
By convention, $a\parallel b$ holds for any pair of element in $X_{0}$.

\begin{definition}
Let $X$ be an omega graph. A \textbf{bridge relation} $R$ on $X$ is a
sequence of ternary relations 
\[
R_{i}\subseteq X_{i}\times X_{i}\times X_{i+1} 
\]
\textbf{\ }

defined for $i\geq 0$ and having the following properties:

\begin{enumerate}
\item  $\left( a,b,c\right) \in R_{i}\Rightarrow a\parallel b$ and $%
dom_{i}^{i+1}c=a$ and $cod_{i}^{i+1}c=b$

\item  $(a,a,id_{i+1}^{i}a)\in R_{i}$ \ for all $i\geq 0$ and for all $a\in
X_{i}$

\item  $(a,b,c),$ $(a,b,d)\in R_{i}\Longrightarrow c=d$
\end{enumerate}
\end{definition}

If $X$ is also an omega magma, then a bridge relation on $X$ is just a
bridge relation on its underlying omega graph $U_{G}^{M}X$. If $X$ is an
omega magma or an omega graph with a bridge relation $R$ then we shall call
the pair $(X,R)$ a \textbf{bridge magma} or a \textbf{bridge graph.}

We note that every omega graph has a diagonal bridge relation denoted by the
symbol $R^{\Delta }$. Here $R_{i}^{\Delta }$ has as its elements all triples
of the form $(a,a,id_{i+1}^{i}a)$.

The significance of this terminology is evident. If $\left( a,b,c\right) \in
R_{i}$ then the $(i+1)$-cell $c$ is a ``bridge'' from the $i$-cell $a$ to
the $i$-cell $b$. One can imagine the parallel cells $a,b$ as being the two
banks of a river. In the literature bridge-cells like $c$ are often part of
what is commonly called a ``contraction'' but we think our terminology is
more descriptive.

\begin{definition}
Let $X$and $Y$be omega magmas and $R^{X},R^{Y}$be bridge relations on $X$
and $Y$. A morphism of omega magmas $F:X\rightarrow Y$ is called a \textbf{%
bridge morphism} if $(a,b,c)\in R_{i}^{X}\Rightarrow (Fa,Fb,Fc)\in R_{i}^{Y}$%
. We shall denote such a bridge morphism by $F:(X,R^{X})\rightarrow
(Y,R^{Y}) $.
\end{definition}

We shall require one more concept before we can define weak omega
categories. This idea is due to Penon and is the essential feature of what
he calls a categorical ``stretching'' [18].

\begin{definition}
A \textbf{categorical Penon morphism} is a bridge morphism between omega
magmas 
\[
F:(X,R^{X})\rightarrow (Y,R^{Y}) 
\]

satisfying the following conditions :

\begin{enumerate}
\item  $Y$ is a strict omega category and $R^{Y}=R^{\Delta }$

\item  $(a,b,c)\in R_{i}^{X}\Longrightarrow Fa=Fb$ and $%
Fc=id_{i+1}^{i}(Fa)=id_{i+1}^{i}(Fb)$

\item  $a,b\in X_{i}$ and $a\parallel b$ and $Fa=Fb\Longrightarrow \exists c$
and $(a,b,c)\in R_{i}^{X}$
\end{enumerate}
\end{definition}

We shall often write $F:X\rightarrow Y$ \ but say that $F$ is a categorical
Penon morphism when the bridge relation $R^{X}$ is understood.

We are now ready to define weak omega categories. The reader will find it
helpful to glance at the following diagram while reading the definition: 
\[
\mathbf{X:\;} 
\begin{array}{ccccc}
& {\small \lambda }^{X} &  & {\small \kappa }^{X} &  \\ 
X_{1} & \leftrightarrows & X_{2} & \rightarrow & X_{3} \\ 
& {\small \rho }^{X} &  &  & 
\end{array}
\]

\begin{definition}
A \textbf{weak omega category} $\mathbf{X}$ consists of the following seven
elements (see the above diagram):

\begin{enumerate}
\item  An omega magma $X_{1}$ called the \textbf{underlying magma} of $%
\mathbf{X}$

\item  An omega magma $X_{2}$ called the \textbf{coherence magma} of $%
\mathbf{X}$

\item  A morphism $\lambda ^{X}:X_{2}\rightarrow X_{1}$ of omega magmas
called the coherence morphism of $\mathbf{X}$

\item  A morphism $\rho ^{X}:X_{1}\rightarrow X_{2}$ of omega \emph{graphs }%
which splits $\lambda ^{X}$ in the category of omega graphs, i.e. 
\[
U_{G}^{M}(\lambda ^{X})\circ \rho ^{X}=1_{U_{G}^{M}X_{1}} 
\]

\item  A bridge relation $R^{X}$ on $X_{2}$

\item  A strict omega category $X_{3}$

\item  A categorical Penon morphism $\kappa ^{X}:(X_{2},R^{X})\rightarrow
(X_{3},R^{\Delta })$
\end{enumerate}
\end{definition}

We next define an omega functor in the obvious way.

\begin{definition}
Let $\mathbf{X}$ and $\mathbf{Y}$ be weak omega categories. An \textbf{omega
functor} $\mathbf{F}:\mathbf{X}\rightarrow \mathbf{Y}$ is a triple of omega
magma morphisms $\mathbf{F}=(F_{1},F_{2},F_{3})$ with the following
properties:

\begin{enumerate}
\item  $F_{i}:X_{i}\rightarrow Y_{i}$

\item  $F_{2}:(X_{2},R^{X})\rightarrow (Y_{2},R^{Y})$ is a morphism of
bridge magmas

\item  The following diagram commutes in the category of omega magmas: 
\[
\begin{array}{cccccc}
&  & {\small \lambda }^{X} &  & {\small \kappa }^{X} &  \\ 
& X_{1} & \leftarrow & X_{2} & \rightarrow & X_{3} \\ 
& \mid &  & \mid &  & \mid \\ 
{\small F}_{1} & \mid & {\small F}_{2} & \mid & {\small F}_{3} & \mid \\ 
& \mid &  & \mid &  & \mid \\ 
& \downarrow & {\small \lambda }^{Y} & \downarrow & \kappa ^{Y} & \downarrow
\\ 
& Y_{1} & \leftarrow & Y_{2} & \rightarrow & Y_{3}
\end{array}
\]

\item  $F_{2}\circ \rho ^{X}=\rho ^{Y}\circ F_{1}$ in the category of omega
graphs.
\end{enumerate}
\end{definition}

\begin{theorem}
The category, Omega\_Cat, in which an object is a weak omega category and a
morphism is an omega functor is the category of models of an essentially
algebraic theory. By theorem 10 it is locally finitely presentable and hence
complete and cocomplete.
\end{theorem}

\noindent \textbf{Proof:}

We only sketch the argument since it is a simple exercise in applying the
definition of an essentially algebraic theory [1].

The sorts of the theory are indexed by pairs of integers $(i,j)$ with $%
i=1,2,3$ and $j\geq 0$. The sort indexed by the pair $(i,j)$ is the sort of $%
j$-cells of the omega magma $X_{i}$.

The total operations are the domain, codomain, and identity functions of
each of the three omega magmas together with the functions defining the
morphisms $\lambda ^{X},\rho ^{X}$ and $\kappa ^{X}$.

The partial operations are the composition operations for each of the three
omega magmas and also the relations $R_{i}^{X}$ comprising the bridge
relation $R^{X}$. The relation $R_{i}^{X}$ is in fact a partial function $%
X_{i}\times X_{i}\rightarrow X_{i+1}$ because it is single valued by
definition 11, number 3.

The axioms of the theory are all equations. These state that the $X_{i}$ are
omega magmas for $i=1,2,3$, that $\lambda ^{X}$ and $\kappa ^{X}$ are
morphisms of omega magmas, that $X_{3}$ is a strict omega category, that $%
R^{X}$ is a bridge relation, that $\kappa ^{X}$ is a categorical Penon
morphism, that $\rho ^{X}$ is a morphism of omega graphs and that $\rho ^{X}$
splits $\lambda ^{X}$ in the category of omega graphs.

It is easy to check that a morphism of models is exactly an omega functor
between weak omega categories. $\blacksquare $

There is an obvious forgetful functor from \textit{Omega\_Cat} to \textit{%
Omega\_Graph} that sends a weak omega category $\mathbf{X}$ to the omega
graph underlying the omega magma $X_{1}$. Both categories are locally
finitely presentable and it is easy to see that this forgetful functor
preserves limits and directed colimits. It follows from theorem 8 that it
has a left adjoint. We have not been able to identify this left adjoint and
suspect it has no simple description. In any case \textit{Omega\_Cat} is
very far from being monadic over \textit{Omega\_Graph}. This is a
consequence of Beck's ``precise tripleability theorem'' [15]. The forgetful
functor from \textit{Omega\_Cat} fails to create split coequalizers (in fact
it creates no colimits whatsoever) and hence by Beck's theorem \textit{%
Omega\_Cat} is not monadic over \textit{Omega\_Graph.} The problem arises
because this forgetful functor forgets far too much structure.

It is now easy to define a weak n-category.

\begin{definition}
A weak omega category $\mathbf{X}$ is a \textbf{weak n-category} if its
underlying omega magma $X_{1}$ is n-skeletal. (Recall that this means that $%
X_{1}$ has only identity cells above dimension n.) \textbf{Weak\_nCat} is
the full subcategory of Omega\_Cat in which an object is a weak n-category.
\end{definition}

The definition makes it obvious that \textit{Weak\_nCat }is itself the
category of models of an essentially algebraic theory and as such is locally
finitely presentable, complete and cocomplete.

\section{Strict omega categories are weak omega categories}

We begin exploring the properties of \textit{Omega\_Cat} by asking whether
strict omega categories as defined in definition 5 are weak omega categories
in the sense of definition 14.

The reader has probably already noticed that this is so, indeed in a trivial
way. Let \textit{Diag\_Omega\_Cat} denote the full subcategory of \textit{%
Omega\_Cat} with objects $\mathbf{X}$ defined by the properties that $%
X_{1}=X_{2}=X_{3}$ , that $\lambda ^{X},\rho ^{X}$ and $\kappa ^{X}$ are all
the identity morphism and that $R^{X}=R^{\Delta }$. Note that for any such
object the omega magmas $X_{i}$ all are equal to the same strict omega
category $X_{3}$. Therefore \textit{Diag\_Omega\_Cat} is isomorphic to 
\textit{Strict\_Cat}, the category of strict omega categories.

The definition of weak omega categories is illuminated in a more interesting
way by considering the full subcategory \textit{Weak\_1Cat} of \textit{%
Omega\_Cat}. Define 
\[
U_{M}^{W}:\text{\textit{Weak\_1\_Cat}}\rightarrow \mathit{Omega\_Magma} 
\]

\noindent to be the functor that sends an object $\mathbf{X}$ to the omega
magma $X_{1}$. Now regard \textit{Cat, }the category of small categories, as
embedded in \textit{Omega\_Magma} by a functor that sends a category $C$ to
the strict omega category $\widehat{C}$ which is identical to $C$ in
dimensions $0$ and $1$ and which is $1$-skeletal. Let $\widehat{\mathit{Cat}}
$ denote this full subcategory of \textit{Omega\_Magma.} Finally, define 
\[
S:\widehat{\mathit{Cat}}\rightarrow \mathit{Omega\_Cat} 
\]

\noindent to be the functor that sends an object of $\widehat{\mathit{Cat}}$
to the obvious object of \textit{Diag\_Omega\_Cat. }This object is clearly
also an object of \textit{Weak\_1Cat.}

\begin{theorem}
The functor $U_{M}^{W}$ has its image in\textit{\ }$\widehat{\mathit{Cat}}$.
In fact it is a retraction split by $S$.
\end{theorem}

\noindent \textbf{Proof:\ }

It will suffice to prove the first statement since the definition of $S$
will then make the second trivial.

We must show that if $\mathbf{X}$ is a weak omega category such that $X_{1}$
is 1-skeletal then the composition law $\bigcirc _{0}^{1}$ and the identity
function $id_{1}^{0}$ for $X_{1}$ define a structure of an ordinary category
with the 0-cells of $X_{1}$ as objects and the 1-cells of $X_{1}$ as
morphisms. We first prove that $\bigcirc _{0}^{1}$ is an associative
composition law.

Let $f,g,h$ be 1-cells of $X_{1}$ such that $f\bigcirc _{0}^{1}g$ and $%
g\bigcirc _{0}^{1}h$ are defined. Since $\rho ^{X}$ is a morphism of omega
graphs we conclude that $\rho ^{X}f\bigcirc _{0}^{1}\rho ^{X}g$ and $\rho
^{X}g\bigcirc _{0}^{1}\rho ^{X}h$ are defined in $X_{2}$. Let $u\equiv (\rho
^{X}f\bigcirc _{0}^{1}\rho ^{X}g)\bigcirc _{0}^{1}\rho ^{X}h$ and $v\equiv
\rho ^{X}f\bigcirc _{0}^{1}(\rho ^{X}g\bigcirc _{0}^{1}\rho ^{X}h)$. Clearly 
$u\parallel v$ and since $X_{3}$ is a strict omega category and $\kappa ^{X}$
is an omega magma morphism we conclude $\kappa ^{X}u=\kappa ^{X}v$. Since $%
\kappa ^{X}$ is a categorical Penon morphism it follows that $\exists c$
such that $(u,v,c)\in R_{1}^{X}$. Now $\lambda ^{X}c$ must be an identity
cell in $X_{1}$ because it is a 2-cell and $X_{1}$ has been assumed to be
1-skeletal. Consequently $\lambda ^{X}u=\lambda ^{X}v$ because $c$ is a
bridge-cell from $u$ to $v$. Since $\lambda ^{X}$ is a morphism of omega
magmas we conclude that 
\[
(\lambda ^{X}\rho ^{X}f\bigcirc _{0}^{1}\lambda ^{X}\rho ^{X}g)\bigcirc
_{0}^{1}\lambda ^{X}\rho ^{X}h=\lambda ^{X}\rho ^{X}f\bigcirc
_{0}^{1}(\lambda ^{X}\rho ^{X}g\bigcirc _{0}^{1}\lambda ^{X}\rho ^{X}h) 
\]

\noindent But $\lambda ^{X}$ is split by $\rho ^{X}$. Thus $(f\bigcirc
_{0}^{1}g)\bigcirc _{0}^{1}h=f\bigcirc _{0}^{1}(g\bigcirc _{0}^{1}h)$ as
desired.

The proofs of the identity laws follow exactly the same pattern as the proof
of the associative law.$\blacksquare \bigskip $

\noindent \textbf{Scholium.} \ We wish to call to the reader's attention the
pattern evident in the preceding proof. This is the standard method for
proving that a desired coherence law must hold in the omega magma $X_{1}$
whenever $\mathbf{X}$ is a weak omega category.

One first lifts the individual $i$-cells which are involved in the coherence
law to $X_{2}$ using the omega graph morphism $\rho ^{X}$. \ It is \textit{%
very important }to note that this is a lifting of the \textit{individual
cells}, \textit{not} of the composites they form in $X_{1}$. One then
reassembles the coherence diagram in $X_{2}$ using the lifts of these cells
from $X_{1}$. Next one observes that certain paths through this diagram
define composite $i$-cells $u$ and $v$ in $X_{2}$ which are parallel.
Applying $\kappa ^{X}$ one infers that the images of $u,v$ are equalized by $%
\kappa ^{X}$ because $X_{3}$ is a strict omega category in which the desired
coherence law holds as an equality. Since $\kappa ^{X}$ is a categorical
Penon morphism one can then deduce that there is a unique $(i+1)$-cell $c$
that is a bridge from $u$ to $v$ in $X_{2}$. The cell $\lambda ^{X}c$ is
then a bridge from $\lambda ^{X}u$ to $\lambda ^{X}v$ in $X_{1}$. The facts
that $\lambda ^{X}$ is an omega magma morphism split by the omega graph
morphism $\rho ^{X}$ then yields the conclusion that the cell $\lambda ^{X}c$
is a bridge between the desired composites in $X_{1}$ in the way illustrated
by the concluding lines of the preceding proof.

In general the bridge-cell $\lambda ^{X}c$ in $X_{1}$ will be neither an
identity nor an isomorphism. However, it will be an equivalence in any
reasonable sense as the following argument shows.

Since $c$ is a bridge from $u$ to $v$ the axioms defining a categorical
Penon morphism assure us that there is also a bridge $c^{\prime }$ from $v$
to $u$. Then both $c\bigcirc _{i}^{i+1}c^{\prime }$ and $c^{\prime }\bigcirc
_{i}^{i+1}c$ are defined and $c\bigcirc _{i}^{i+1}c^{\prime }\parallel
id_{i+1}^{i}\left( u\right) $ while $c^{\prime }\bigcirc
_{i}^{i+1}c\parallel id_{i+1}^{i}\left( v\right) $. One again appeals to the
fact that $\kappa ^{X}$ is a categorical Penon morphism to deduce that $%
\kappa ^{X}$ equalizes both pairs of parallel cells. Consequently there is a
unique $(i+2)$-cell $d$ that is the bridge from $c\bigcirc
_{i}^{i+1}c^{\prime }$ to $id_{i+1}^{i}(u)$. Applying $\lambda ^{X}$ to
these cells we find that $\lambda ^{X}d$ is a bridge from $\lambda
^{X}c\bigcirc _{i}^{i+1}\lambda ^{X}c^{\prime }$ to $id_{i+1}^{i}(\lambda
^{X}u)$. Thus, while the $(i+1)$-cell $\lambda ^{X}c$ is not an isomorphism
there is a cell $\lambda ^{X}c^{\prime }$ such that the composite $\lambda
^{X}c\bigcirc _{i}^{i+1}\lambda ^{X}c^{\prime }$ and the identity cell $%
id_{i+1}^{i}(\lambda ^{X}u)$ are bridged by an $(i+2)$-cell$\ \lambda ^{X}d$%
. Of course, the $(i+2)$-cell$\ \lambda ^{X}d$ is itself neither an identity
nor an isomorphism. It does however have a ``quasi-inverse'' as did $\lambda
^{X}c$. The composite of $\lambda ^{X}d$ with its quasi-inverse can then be
connected to $id_{i+2}^{i+1}(\lambda ^{X}u)$ by yet another bridge cell of
dimension $(i+3)$ and so forth. $\blacksquare $

\section{Bicategories and weak 2-categories}

For a definition of bicategories and a discussion of some of their
properties the reader may consult [6,7,12,14,15,16]. Throughout this section
we shall maintain as much consistency as possible with the bicategorical
notation of [12].

Let \textit{Bicat} denote the category in which an object is a bicategory
and a morphism is strong homomorphism between bicategories (i.e. a morphism
that preserves the operations and identities ''on the nose'', not just up to
isomorphism; in [12] this is called a strict homomorphism). Our aim in this
section is to show that \textit{Bicat} is a retract of \textit{Weak\_2\_Cat}%
. Thus we shall construct a pair of functors 
\[
\begin{array}{ccc}
& {\small U}_{B}^{W} &  \\ 
\text{\textit{Weak\_2\_Cat}} & \rightleftarrows & \text{\textit{Bicat}} \\ 
& {\small S}_{W}^{B} & 
\end{array}
\]

\noindent with the property that $U_{B}^{W}\circ S_{W}^{B}\backsimeq 1_{%
\text{\textit{Bicat}}}$.

We first construct the functor $U_{B}^{W}$.

\begin{theorem}
Let $\mathbf{X}$ be a weak 2-category. Then its underlying omega magma, $%
X_{1}$, can be structured as a bicategory for which the coherence data
consists of 2-cells which are the images under $\lambda ^{X}$ of
bridge-cells in $X_{2}$.
\end{theorem}

\noindent \textbf{Proof:}

In Section 7.1 we shall define for any weak omega category $\mathbf{X}$ and
any pair of $(i-1)$-cells $a,b\in X_{0}$ a weak omega category $\mathbf{X(}%
a,b)$ in which a 0-cell $f$ is an $i$-cell of $X_{1}$ with domain $a$ and
codomain $b$. The composition laws of $\mathbf{X(}a,b)$ are the suitably
reindexed restrictions of the composition laws of $\mathbf{X}$. However,
when analyzing ``hom categories'' like $\mathbf{X(}a,b)$ it is usually
convenient to retain the indexing of the composition laws of $\mathbf{X}$
when discussing composition in $\mathbf{X}(a,b)$.

In the situation at hand, for any pair of 0-cells $a,b\in (X_{1})_{0}$ we
can consider the weak omega category $\mathbf{X}(a,b)$. Since $\mathbf{X}$
is 2-skeletal we know that $\mathbf{X}(a,b)$ is 1-skeletal. Theorem 18 then
asserts that $X(a,b)_{1}$ is in fact an ordinary category whose composition
law is the restriction of the law $\bigcirc _{1}^{2}$ in $X_{1}$. Thus the
first requirement $X_{1}$ must satisfy as a bicategory is met.

Our next task is to identify, for each triple of 0-cells $a,b,c\in
(X_{1})_{0}$ the composition functor 
\[
C_{abc}:X(a,b)_{1}\times X(b,c)_{1}\rightarrow X(a,c)_{1} 
\]

\noindent If $(f,g)$ is an ordered pair of 0-cells (i.e. a pair of 1-cells
of $X_{1}$) we define $C_{abc}(f,g)=f\bigcirc _{0}^{1}g$ where the operation
is that of $X_{1}$. Similarly, if $(\beta ,\gamma )$ is a pair of 1-cells
(i.e. a pair of 2-cells of $X_{1}$) we define $C_{abc}(\alpha ,\beta
)=\alpha \bigcirc _{0}^{2}\beta $ where the operation is again that of $%
X_{1} $.

We next show that $C_{abc}$ defined in this way is in fact a functor. This
is equivalent to showing that the interchange law holds as an equation in $%
X_{1}$, i.e. that for a 4-tuple $(\beta ,\gamma ,\beta ^{\prime },\gamma
^{\prime })$ of compatible 2-cells in $X_{1}$we have the equation (in $X_{1}$%
) 
\[
(\beta \bigcirc _{0}^{2}\gamma )\bigcirc _{1}^{2}(\beta ^{\prime }\bigcirc
_{0}^{2}\gamma ^{\prime })=(\beta \bigcirc _{1}^{2}\beta ^{\prime })\bigcirc
_{0}^{2}(\gamma \bigcirc _{1}^{2}\gamma ^{\prime }) 
\]

\noindent It will suffice to exhibit a 3-cell of $X_{1}$ that is a
bridge-cell from the left hand side of this equation to its right hand side.
\noindent For by hypothesis, $X_{1}$ is 2-skeletal and hence any such
bridge-cell must be an identity cell.

We produce such a bridge-cell by following the method described in the
Scholium in Section 3. First use the omega graph morphism $\rho ^{X}$ to
lift these four 2-cells from $X_{1}$ to $X_{2}$. Then in $X_{2}$ assemble
the two composites which correspond to those appearing in the last equation.
Notice that these are parallel 2-cells in $X_{2}$ and that since $\kappa
^{X} $ is an omega magma morphism and $X_{3}$ is a strict omega category
this pair of composite 2-cells is equalized by $\kappa ^{X}$. It follows
that this pair of composite cells is bridged by a 3-cell $\sigma $ in $X_{2}$%
. Next note that since $\rho ^{X}$ splits $\lambda ^{X}$ and since the
latter is an omega magma morphism the images under $\lambda ^{X}$ of the two
composite 2-cells in $X_{2}$ are precisely the 2-cells that appear on the
right and left hand sides of the last equation. Consequently the 3-cell $%
\lambda ^{X}\sigma $ is the desired bridge cell and must be an identity cell
because $X_{1}$ is 2-skeletal.

Our next task is to produce the coherence data that makes $X_{1}$ a
bicategory.

For any 4-tuple $(a,b,c,d)$ of 0-cells of $X_{1}$ there are two obvious
functors: 
\[
\begin{array}{ccc}
& {\small C}_{abd}{\small \circ (1\times C}_{bcd}{\small )} &  \\ 
X(a,b)_{1}\times X(b,c)_{1}\times X(c,d)_{1} & \rightrightarrows & X(a,d)_{1}
\\ 
& {\small C}_{acd}{\small \circ (C}_{abc}{\small \times 1)} & 
\end{array}
\]

\noindent We must construct a natural isomorphism 
\[
\alpha :{\normalsize C}_{abd}{\normalsize \circ (1\times C}_{bcd}%
{\normalsize )\rightarrow C}_{acd}{\normalsize \circ (C}_{abc}{\normalsize %
\times 1)} 
\]
called the associator of the bicategory. For any triple of 0-cells $%
(f,g,h)\in $ $X(a,b)_{1}\times X(b,c)_{1}\times X(c,d)_{1}$ (these are
1-cell in $X_{1}$) we must define a 1-cell isomorphism (a 2-cell in $X_{1}$) 
\[
\alpha _{fgh}:f\bigcirc _{0}^{1}(g\bigcirc _{0}^{1}h)\rightarrow (f\bigcirc
_{0}^{1}g)\bigcirc _{0}^{1}h) 
\]
(where the operations are the operations of $X_{1}$). This is easily done
using the technique of the Scholium\textbf{\ }of Section 3. One first uses
this technique to exhibit a bridge-cell between these two composite 1-cells
in $X_{1}$. This bridge-cell must be an isomorphism because its composite
with the bridge cell in the opposite direction is connected to an 2-cell
identity by a three-cell. This three cell must itself be an identity cell
because $X_{1}$ is 2-skeletal. Consequently the original bridge-cell is an
isomorphism.

That these 2-cell isomorphisms are the components of a natural isomorphism
of functors is also easy to prove using the method described in the
Scholium. One simply lifts the 2-cells of $X_{1}$ comprising the diagrams
which we must show are commutative to $X_{2}$. After reassembling these
diagrams in $X_{2}$ the by-now-standard argument shows that the two paths
through this diagram define composite 2-cells in $X_{2}$ that are bridged by
a 3-cell. Since $X_{1}$ is 2-skeletal the image of this 3-cell bridge under $%
\lambda ^{X}$ is an identity 3-cell in $X_{1}$ and so the diagram must
commute in $X_{1}$.

Precisely the same arguments produce natural isomorphisms 
\[
l_{ab}:C_{aab}\circ (I_{a}x1)\rightarrow 1_{X(a,b)_{1}} 
\]
\[
r_{ab}:C_{abb}\circ (1\times I_{a})\rightarrow 1_{X(a,b)_{1}} 
\]

\noindent In these last two equations the functors $I_{a}\times 1$ and $%
1\times I_{a}$ are respectively left and right composition of 1-cells in $%
X_{1}$ (0-cells in $X(a,b)_{1}$) with the identity 1-cell $id_{1}^{0}(a)$. $%
\blacksquare $

It is easy to see that when $\mathbf{X}$ and $\mathbf{Y}$ are weak
2-categories an omega functor $\mathbf{F}:\mathbf{X}\rightarrow \mathbf{Y}$
induces a strong homomorphism of the bicategory structures we have just
constructed on $X_{1}$ and $Y_{1}$. This follows from the fact that the
omega magma morphism $F_{2}$ is required to be a morphism of bridge magmas.

This completes the construction of the functor $U_{B}^{W}$.

We next construct the functor $S_{W}^{B}$ from \textit{Bicat} to \textit{%
Weak\_2\_Cat}.

We begin by noting that any bicategory can be regarded as a 2-skeletal omega
magma paired with associated coherence data. Forgetting this coherence data
as well as the composition laws then defines a forgetful functor $U_{G}^{B}$
from \textit{Bicat} to \textit{Omega\_Graph}. Now \textit{Bicat} is clearly
the category of models of an essentially algebraic theory and it is easy to
see that this forgetful functor creates limits and directed colimits and
hence preserves them. \ Then theorem 8 tells us that $U_{G}^{B}$ has a left
adjoint $L_{B}^{G}$.

Now let $Y$ be a bicategory regarded as a 2-skeletal omega magma provided
with coherence data. We shall define a diagram of omega magmas, $%
S_{W}^{B}Y\equiv \mathbf{X}$, and show that this diagram is in fact a weak
2-category.

Let $X_{1}\equiv Y$. Let $X_{2}\equiv L_{B}^{G}\circ U_{G}^{B}(Y)$, the
bicategory freely generated by the omega graph underlying $Y$. Define $\rho
^{X}\equiv \eta U_{G}^{B}$, the component of the unit of the $%
L_{B}^{G}\dashv U_{G}^{B}$ adjunction corresponding to $Y$. Define $\lambda
^{X}$ $\equiv $ $\varepsilon _{Y}$, the component of the counit of this
adjunction corresponding to $Y$. Then $\rho ^{X}$ splits $\lambda ^{X}$ in
the category of omega graphs by the triangle identities for the adjunction.

Next define $X_{3}\equiv L_{S}^{G}\circ $ $U_{G}^{B}(Y)$, the strict omega
category freely generated by the omega graph $U_{G}^{B}Y$. Clearly this
strict omega category is 2-skeletal and is therefore a freely generated,
strict 2-category. Define $\kappa ^{X}\equiv \eta U_{M}^{S}$, the component
of the unit for the $L_{S}^{M}\dashv U_{M}^{S}$ adjunction between \textit{%
Strict\_Cat} and \textit{Omega\_Magma} that corresponds to the omega magma
underlying $X_{2}\equiv L_{B}^{G}\circ U_{G}^{B}(Y)$. In less mysterious
terms $\kappa ^{X}$ is an omega magma morphism from the free bicategory to
the free strict 2-category generated by the same omega graph. This morphism
is induced by the omega magma congruence generated by the relations that
must hold in every strict 2-category.

It should be clear that $X$ depends functorially on the bicategory $Y$. Thus
we can complete the construction of $S_{C}^{B}$ by proving:

\begin{theorem}
$\mathbf{X}$ is a weak 2-category and the coherence data for $Y$ is
constructed from the images of the bridge-cells of $X_{2}$ under $\lambda
^{X}$.
\end{theorem}

\noindent \textbf{Proof}:

We shall first construct a bridge relation $R^{X}$ on $X_{2}$, the
bicategory freely generated by the (2-skeletal) omega graph underlying $Y$.
This construction is made possible by the coherence theorem for bicategories
[12,14,16]. It will be evident from the construction that the 2-cell
coherence isomorphisms for $Y$ are the images under $\lambda ^{X}$ of bridge
cells defined by $R^{X}$.

For $i\neq 1$ define $R_{i}^{X}\equiv \left\{ (a,a,id_{i+1}^{i}a)|a\in
(X_{2})_{i}\right\} $.

For $i=1$ define 
\[
R_{1}^{X}\equiv \left\{ (a,b,c)|a,b\in (X_{2})_{1},c\in (X_{2})_{2}\text{ , }%
c\text{ an isomorphism from }a\text{ to }b\right\} 
\]
The fact that $X_{2}$ is freely generated coupled with the coherence theorem
for bicategories guarantees that 
\[
\left( a,b,c\right) ,\left( a,b,d\right) \in R_{1}^{X}\Rightarrow c=d 
\]

\noindent because the only 2-cell isomorphisms in $(X_{2})_{2}$ are
composites of coherence isomorphisms. (The point here is that freeness
guarantees that the only automorphism of a 1-cell is the identity
automorphism.) Moreover, such isomorphisms must necessarily connect parallel
1-cells. It follows that $R^{X}$ $\ $is a bridge relation. Since $\lambda
^{X}$ is defined as a counit in \textit{Bicat} it a strong homorphism of
bicategories and hence necessarily sends the bridge elements determined by $%
R^{X}$ to coherence isomorphisms of $X_{1}=Y$.

It remains to prove that $\kappa ^{X}$ is a categorical Penon morphism for
this bridge relation. If $(a,b,c)\in R_{1}^{X}$ then $a,b$ are obviously
equalized by $\kappa ^{X}$. Moreover, the fact that $X_{3}$ is freely
generated means that there are no relations among its 2-cells so that $%
\kappa ^{X}$ must send any 2-cell isomorphism to an identity cell. On the
other hand the coherence theorem for bicategories tells us that if $%
a\parallel b$ are 1-cells of $X_{2}$ which are equalized by $\kappa ^{X}$
then there is a unique 2-cell isomorphism in $X_{2}$ from $a$ to $b$. The
fact that $X_{2}$ is freely generated means that if $c$ is this isomorphism
then $\left( a,b,c\right) \in $ $\kappa ^{X}$.

To verify the Penon condition for $R_{i}^{X}$, $i\neq 1$ we observe that by
construction $\kappa ^{X}$ is an isomorphism on $i$-cells for $i=0$ and for $%
i\geq 3$. \ We can thus complete the proof by showing that $\kappa ^{X}$ is
faithful when restricted to 2-cells, i.e. that it does not equalize any
parallel pair of 2-cells in $X_{2}$. This fact will follow from the
following easy lemma.

We first recall some terminology. A \textit{clique} is an ordinary category
equivalent to the terminal category; in other words, every pair of objects
in a clique is connected by a unique isomorphism.

\begin{lemma}
Let $C$ be an ordinary category for which $ISO(C)$, the maximal subgroupoid
of $C$ whose morphisms are all isomorphisms,\ is a coproduct (disjoint
union) of cliques. Let $\widetilde{C}$ denote the category obtained from $C$
by identifying isomorphic objects. Then the projection functor $C\rightarrow 
\widetilde{C}$ is fully faithful.
\end{lemma}

\noindent \textbf{Proof of Lemma:}

The hypothesis ensures that $\widetilde{C}$ is isomorphic to the category
whose objects are the isomorphism classes $[a]$ of $C$ and in which $%
\widetilde{C}([a],[b])\equiv C(a,b)$, the latter sets being canonically
isomorphic for any choice of representatives by the clique hypothesis. $%
\blacksquare $

We apply the preceding lemma to the category $X_{2}(a,b),a,b$ 0-cells of $%
X_{2}$. The fact that $X_{2}$ is freely generated as a bicategory coupled
with the coherence isomorphism for bicategories ensures that the hypothesis
of the lemma is satisfied. We conclude that $\kappa ^{X}$ is faithful when
restricted to 2-cells thus completing the proof of the theorem. $%
\blacksquare $

\section{The weak omega categories of Penon and Batanin}

In this section we shall show that the weak omega categories defined by
Penon in [18] are instances of our weak omega categories. \ We shall also
offer an informal argument to support our claim that this is also true for
the weak omega categories defined by Batanin [5]. In our opinion only our
own lack of familiarity with the technical details of Batanin's work
prevents us from making this argument rigorous.

\subsection{The category \textit{Prolixe}}

Let us begin by considering the category \textit{Pen\_Mor} in which an
object is a categorical Penon morphism (definition 13) and in which a
morphism $(H_{W},H_{Z})$ is a commutative square: 
\[
\begin{array}{ccccc}
&  & F &  &  \\ 
& (X,R^{X}) & \rightarrow & (Y,R^{\Delta }) &  \\ 
H_{W} & \downarrow &  & \downarrow & H_{Z} \\ 
& (W,R^{W}) & \rightarrow & (Z,R^{\Delta }) &  \\ 
&  & G &  & 
\end{array}
\]

\noindent It is easy to see that \textit{Pen\_Mor} is the category of models
of an essentially algebraic theory and as such is locally finitely
presentable, complete and cocomplete.

There is a forgetful functor $U_{G}^{P}:$\textit{Pen\_Mor}$\rightarrow $%
\textit{Omega\_Graph} which sends an Penon morphism (an object) to the omega
graph which underlies the omega magma of it domain. This functor clearly
preserves limits and directed colimits (in fact all filtered colimits) and
so has a left adjoint $L_{P}^{G}$. Let $\Bbb{T\equiv }U_{G}^{P}\circ
L_{P}^{G}$ denote the resulting monad on \textit{Omega\_Graph.} Penon
defines a \textbf{prolixe} to be a $\Bbb{T}$ algebra and a morphism of
prolixes to be a $\Bbb{T}$ algebra morphism. There is therefore a category 
\textit{Prolixe} and this is the category of Penon's weak omega categories.
Since $\Bbb{T}$ preserves filtered colimits \textit{Prolixe} is cocomplete.
However, it is not complete since it does not contain Cartesian products.
Thus \textit{Prolixe} cannot be the category of models of any essentially
algebraic theory.

Let us examine the left adjoint $L_{P}^{G}$ in a little more detail. If $Y$
is an omega graph $L_{P}^{G}Y$ is a morphism of bridge magmas 
\[
L_{P}^{G}Y:\text{domain }L_{P}^{G}Y\rightarrow \text{codomain }L_{P}^{G}Y 
\]

\noindent Using the universal property of the left adjoint one sees that
codomain $L_{P}^{G}Y$ is just the strict omega category freely generated by
the omega graph $Y$.

Here is an informal construction of the bridge magma which is the domain of
the categorical Penon morphism $L_{P}^{G}Y$.\ (See part 2 of [18] for a
formal construction.). Let $M(Y)_{0}$ denote the omega magma freely
generated by $Y$ and let $C(Y)$ denote the strict omega category freely
generated by $Y$. Then $\eta _{0}:M(Y)_{0}\rightarrow C(Y)$ denotes the
canonical morphism (unit of the omega magma-strict omega category
adjunction). Begin the construction by formally adjoining a single 2-cell to 
$M(Y)_{0}$ for each ordered pair of parallel one cells equalized by $\eta
_{0}$. Extend $\eta _{0}$ over the new 2-cells by mapping them to the
appropriate identities in $C(Y)$. Define $R_{1}$ to consist of those triples
whose first two coordinates is an ordered pair of equalized, parallel
1-cells and whose third coordinate is the corresponding new 2-cell. Next
``magmify'' the resulting omega graph by freely adjoining all compositions
with the new 2-cells and appropriate identity cells and extend $\eta _{0}$
in the obvious way over this new omega magma which we shall call $M(Y)_{1}$
(note that certain relations must also be added that say that composites
consisting only of previously existing cells remain unaltered). \ Note that $%
M(Y)_{0}\subset M(Y)_{1}$ as omega magmas and that these two omega magmas
have the same 0- and 1-cells. We thus have a morphism $\eta
_{1}:M(Y)_{1}\rightarrow C(Y)$. We repeat this process and obtain a new
diagram $\eta _{2}:M(Y)_{2}\rightarrow C(Y)$ in which $M(Y)_{2}$ has the
same 0-, 1- and 2-cells as $M(Y)_{1}$, $M(Y)_{1}\subset M(Y)_{2}$ and $R_{2}$
has been defined. Continuing in this way and letting $M_{P}(Y)$ denote the
colimit of the obvious directed diagram we obtain domain $L_{P}^{G}Y\equiv
M_{P}(Y)$ and the morphism $L_{P}^{G}Y\equiv \eta
_{P}(Y):M_{P}(Y)\rightarrow C(Y)$.

\subsection{\textit{Prolixe} is a retract}

We shall now show that \textit{Prolixe} is a retract of full subcategory of 
\textit{Omega\_Cat}.

Define \textit{Pen\_Omega\_Cat} to be the full subcategory of \textit{%
Omega\_Cat} in which an object $\mathbf{X}$ has the following properties:

\begin{enumerate}
\item  Define $Y\equiv U_{G}^{M}X_{1}$ then $\kappa ^{X}=\eta _{P}(Y)\equiv
L_{P}^{G}Y$

\item  $\rho ^{X}$ is the unit of the \textit{Omega\_Graph }- \textit{%
Pen\_Mor} adjunction.
\end{enumerate}

There is a functor 
\[
U_{PRO}^{POC}:\text{\textit{Pen\_Omega\_Cat}}\rightarrow \text{\textit{%
Prolixe}} 
\]

\noindent that sends an object $\mathbf{X}$ to the $\Bbb{T}$ algebra $%
U_{G}^{M}\lambda ^{X}:U_{G}^{M}M_{P}(Y)\rightarrow U_{G}^{M}X_{1}$. The
definition of $U_{PRO}^{POC}$ on morphisms is the obvious one.

In the opposite direction, define a functor 
\[
S_{POC}^{PRO}:\text{\textit{Prolixe}}\rightarrow \text{\textit{%
Pen\_Omega\_Cat}} 
\]

\noindent as follows. \ If $h:\Bbb{T}Y$ $\rightarrow Y$ is a $\Bbb{T}$
algebra then $S_{POC}^{PRO}(h)$ is the following weak omega category $%
\mathbf{X}$:

\begin{enumerate}
\item  The unit of $\Bbb{T}$ gives $Y$ the structure of an omega magma
because this unit necessarily splits $h$. Define $X_{1}\equiv Y$

\item  Define $X_{2}\equiv $ $M_{P}(Y)$ and $\rho ^{X}$ to be the unit of $%
\Bbb{T}$ and $\lambda ^{X}\equiv h$. Note that this makes sense because $h $
is necessarily an omega magma morphism by the definition of the omega magma
structure on $Y$.

\item  Define $\kappa ^{X}\equiv \eta _{P}(Y):M_{P}(Y)\rightarrow C(Y)$.
\end{enumerate}

Clearly $U_{PRO}^{POC}\circ S_{POC}^{PRO}=1_{\Pr \text{olixe}}$ and this
shows that \textit{Prolixe} is a retract of \textit{Pen\_Omega\_Cat.}

\subsection{Remarks on Batanin's weak omega categories}

It is our view that it should be possible to adapt the preceding
construction to show that the category of Batanin's weak omega categories
[5] is a retract of some suitable subcategory of \textit{Omega\_Cat} . Let 
\textit{Batanin\_Cat} denote the category of algebras for Batanin's initial,
contractible, higher dimensional operad with a system of contractions. This
category is equivalent to the category of algebras for the monad $\Bbb{B}$
canonically associated with the operad. The unit of the monad gives each
algebra the structure of an omega graph. \ The system of compositions gives
each algebra the structure of an omega magma. Thus the algebra $h:\Bbb{B}%
Y\rightarrow Y$ should correspond as above to the coherence arrow $\lambda
^{X}$ of a weak omega category. Now the category of strict omega categories
is the category of algebras for Batanin's \textit{terminal} operad and for
its associated monad $\Bbb{S}$. The fact that $\Bbb{S}$ is defined by the
terminal operad should give rise to a canonical morphism of omega magmas $%
\Bbb{B}Y\rightarrow \Bbb{S}Y$ . \ The omega magma $\Bbb{S}Y$ is a strict
omega category and this last morphism should be a categorical Penon morphism
because the operad defining $\Bbb{B}$ is contractible. This construction
should define a functor from \textit{Batanin\_Cat} to \textit{Omega\_Cat}.
The functor in the opposite direction should be a simply defined version of
the forgetful functor as it was in the case of \textit{Prolixe.}

\section{The stabilization conjecture}

In [4] Baez and Dolan informally discuss a number of desirable properties a
good theory of weak omega categories might have. Among these is a certain
stability property which we now explain.

\begin{definition}
A weak $(n+k)$ category $\mathbf{X}$ is called $k$\textbf{-tuply monoidal}
if each of the omega magmas $X_{1},X_{2},X_{3}$ has exactly one cell in each
dimension $\leq k-1$. We allow the possibility $n=\omega $. Let $\mathbf{nCAT%
}_{k}$ denote the full subcategory of Omega\_Cat in which an object is a $k$%
-tuply monoidal weak $(n+k)$ category.
\end{definition}

Baez and Dolan suggest that in a good theory of weak omega categories one
should be able to construct a weak omega category \textit{nBD}$_{k}$whose
0-cells are k-tuply monoidal objects like the objects of our category 
\textit{nCat}$_{k}$. They hope that there would then be an omega functor 
\[
S:\mathit{nBD}_{k}\rightarrow \mathit{nBD}_{k+1} 
\]
(the stabilization functor) that would be some sort of equivalence provided $%
k\geq 2$.

We prove a stronger version of this stabilization property here.

\begin{theorem}
For all $n\geq 0$ and $k\geq 2$ there is a functor $S_{n,k}:nCAT_{k}$ $%
\rightarrow nCAT_{k+1}$ that is an isomorphism of categories.
\end{theorem}

\noindent \textbf{Proof:}

The theorem will follow immediately from a similar result in the category of
omega magmas.

A k-tuply monoidal omega magma is an omega magma that has exactly one cell
in each dimension $\leq k-1$. \ Let $M$ be a k-tuply monoidal omega magma
and assume $k\geq 2$. Define $W(M)$ to be the omega graph with a unique
0-cell $\ast $ and all of its higher dimensional identities and whose other $%
i$-cells are the $(i-1)$-cells of $M$ with the obvious reindexing of the
domain, codomain and identity functions of $M$.

Give $W(M)$ the structure of an omega magma in the following way. If $a,b$
are $j$-cells which are compatible for composition over $i$-cells in $W(M)$
then $a,b$ are by definition $(j-1)$-cells of $M$ which are compatible for
composition over $(i-1)$-cells in $M$. Then define the partial composition $%
\bigcirc _{i}^{j}$in $W(M)$ by setting $\bigcirc _{i}^{j}\equiv \bigcirc
_{i-1}^{j-1}$ where the latter law is composition in $M$. This defines
composition laws $\bigcirc _{i}^{j}$ in $W(M)$ for $i\geq 1$. Define $%
\bigcirc _{0}^{j}\equiv \bigcirc _{1}^{j}$ where the right hand side is the
law already defined in $W(M)$. Finally, define $id_{j}^{0}(\ast )\bigcirc
_{0}^{j}a\equiv a$ for all $j$-cells $a$ and similarly for right composition
with this identity.

It is easy to see that $W(M)$ is an omega magma and is just $M$ ``shifted up
one dimension''. This $W$ construction has an obvious inverse if $k\geq 3$
defined by ``forgetting the unique 0-cell $\ \ast $ and all of its
identities''. Clearly both the $W$ construction and its inverse are
functorial.

The theorem now follows from the observation that the $W$ construction can
be applied term by term to any k-tuply monoidal weak $(n+k)$ category and to
the morphisms defining the category and to the bridge relation. The fact
that the resulting functor $S_{n,k}$ is an isomorphism follows from the
existence of the functorial inverse to the $W$ construction for $k\geq 3$ $%
\blacksquare $

\section{Constructing weak omega categories}

In this section we shall discuss some methods for constructing weak omega
categories from other weak omega categories, from omega magmas with special
properties and from omega graphs.

The most basic observation is that Omega\_Cat is the category of models of
an essentially algebraic theory and is locally finitely presentable. \ It is
therefore complete and cocomplete with respect to ordinary conical limits
and colimits.

\subsection{The weak omega category $\mathbf{X}(a,b)$}

Next we construct a weak omega category $\mathbf{X}(a,b)$ for any weak omega
category $\mathbf{X}$ and any pair of $(i-1)$-cells $a,b$ of $X_{1}$. The
0-cells of the underlying omega magma $X(a,b)_{1}$ of $\ \mathbf{X}(a,b)$
will be the $i$-cells of $X_{1}$ with domain $a$ and codomain $b$.

\begin{enumerate}
\item  Define the omega magma $X(a,b)_{1}$ by defining the $j$-cells of $%
X(a,b)_{1}$ to be the $(j+i)$-cells of $X_{1}$ which are mapped by $%
dom_{i-1}^{j+i}$ to $a$ and by $cod_{i-1}^{j+i}$ to $b$. The domain,
codomain and identity functions of $X(a,b)_{1}$ are just those of $X_{1}$
restricted and reindexed in the obvious way. The partial composition
operations of $X(a,b)_{1}$ are again those of $X_{1}$ restricted and
reindexed in the obvious way.

\item  Define the omega magma $X(a,b)_{2}$ by defining the $j$-cells of $%
X(a,b)_{2}$ to be the $(j+i)$-cells of $X_{2}$ which are mapped by $%
dom_{i-1}^{j+i}\circ \lambda ^{X}$ to $a$ and by $cod_{i-1}^{j+i}\circ
\lambda ^{X}$ to $b$. The remaining structure on $X(a,b)_{2}$ is defined in
just the same way as it was for $X(a,b)_{1}$.

\item  Clearly the omega magma morphism $\lambda ^{X(a,b)}$ and the omega
graph morphism $\rho ^{X(a,b)}$ can be defined as the restrictions of the
morphisms $\lambda ^{X}$ and $\rho ^{X}$ and will have the required
properties.

\item  Define the strict omega category $X(a,b)_{3}$ to have as its $j$%
-cells \textit{all} of the $(j+i)$-cells of $X_{3}$. The remaining structure
can be defined as it was for $X(a,b)_{1}$

\item  The bridge relation $R^{X(a,b)}$ on $X(a,b)_{2}$ is just the bridge
relation $R^{X}$ restricted and reindexed in the obvious way.

\item  Finally, the omega magma morphism $\kappa ^{X(a,b)}$ is just the
obvious restriction of $\kappa ^{X}$ and is clearly a categorical Penon
morphism.
\end{enumerate}

This construction of the weak omega category $\mathbf{X}(a,b)$ is functorial
in the following sense. Let $\mathbf{F:X\rightarrow Y}$ be an omega functor
and $a,b$ a pair of $(i-1)$-cells of $X_{1}$. Then $\mathbf{F}$
``restricts'' in the obvious way to an omega functor 
\[
\mathbf{F}(a,b):\mathbf{X}(a,b)\rightarrow \mathbf{Y}(F_{1}a,F_{1}b) 
\]

\subsection{Weak omega categories from omega graphs}

\subsubsection{The functor $P$}

It can be useful to have a method for associating in a functorial way weak
omega categories with omega graphs and with certain diagrams of omega graphs.

In Section 5 we constructed a functor $L_{P}^{G}$, left adjoint to a
forgetful functor, from the category of omega graphs to the category of
categorical Penon morphisms. If $Y$ is an omega graph recall that we denoted 
$L_{P}^{G}Y$ by 
\[
L_{P}^{G}Y\equiv \eta _{P}(Y):M_{P}(Y)\rightarrow C(Y) 
\]
We shall now show how this functor can be used to construct a functor from
omega graphs to weak omega categories.

We shall define a functor 
\[
P:\mathit{Omega\_Graph}\rightarrow \mathit{Omega\_Cat} 
\]

\noindent and call the weak omega category $P(Y)$ the Penon category
associated with $Y$. Temporarily denote $P(Y)$ by $\mathbf{X}$. Then define:

\begin{enumerate}
\item  $X_{1}=X_{2}\equiv M_{P}(Y),\;\lambda ^{X}=\rho ^{X}\equiv
1_{M_{P}(Y)} $

\item  $X_{3}\equiv C(y),\;\kappa ^{X}\equiv \eta _{P}(Y)$
\end{enumerate}

\noindent It should be obvious that this definition extends to morphisms
without difficulty. As far as we have been able to determine $P$ is not left
adjoint to any functor resembling a forgetful functor defined on all of 
\textit{Omega\_Cat}.

\subsubsection{The functor $P_{A}^{\rightarrow }$}

There is a variant of the Penon category construction which we shall find
useful in Section 9. Recall the category of categorical Penon morphisms, 
\textit{Pen\_Mor }from the beginning of Section 5\textit{.} An object is a
categorical Penon morphism (definition 13) and a morphism is the obvious
commutative square. Fix a strict omega category $A$ and consider the
subcategory \textit{Pen\_Mor(A)} in which every object has codomain $A$ and
in which every morphism is the identity functor of $A$ on its codomain leg.
Consider the subcategory of \textit{Omega\_Graph}$^{\rightarrow }$, the
arrow category of \textit{Omega\_Graph,} in which every object has codomain $%
U_{G}^{C}A$ , the omega graph underlying $A$ and in which every morphism has
the identity of $U_{G}^{C}A$ as its codomain leg. Denote this subcategory by 
\textit{Omega\_Graph}$^{\rightarrow UA}$. There is an obvious forgetful
functor 
\[
U_{GMA}^{PMA}:\text{\textit{Pen\_Mor(A)}}\rightarrow \text{\textit{%
Omega\_Graph}}^{\rightarrow UA} 
\]

\noindent We shall argue that this functor has a left adjoint. From an
informal point of view this is obvious from the heuristic construction of
the omega magma morphism $\eta _{P}:$\ $M(Y)_{P}\rightarrow C(Y)\,\ $\ in
Section 5. This construction started with the canonical arrow $\eta _{0}$
from the omega magma freely generated by the omega graph $Y$ to the strict
omega category generated by $Y$. But the construction would clearly make
sense starting with any arrow from this free omega magma to any strict omega
category; this is the data provided by objects of the category \textit{%
Omega\_Graph}$^{\rightarrow UA}$.

A formal argument for the existence of the left adjoint runs as follows. The
forgetful functor obviously creates limits and directed colimits and hence
preserves them. (One should note that the categorical product in each
category is in fact the pullback over the codomain.) To apply theorem 8 we
would like to show that both categories in question are locally finitely
presentable. This can be done in two different ways.

One can show that each of these categories is the category of models of an
essentially algebraic theory. This is rather tedious because one must index
sorts in an unusual way. Sorts are types of cells in the domain magma or
graph. Such a sort is indexed by its dimension and also by its image cell in
the codomain. This leads to complicated bookkeeping when it comes to
defining total and partial operations.

A more straightforward approach is to apply theorem 9 which characterizes
locally finitely presentable categories. It is easy to show that each of
these categories is cocomplete. If we can show that each has a strong
generator then it will follow that each is locally finitely presentable. The
idea now is to use the strong generators of the ambient categories (which
exist because the ambient categories are locally finitely presentable) to
construct the desired strong generators.\ We illustrate this construction
for \textit{Pen\_Mor(A)}.

Let \textsf{G} be strong generator of \textit{Pen\_Mor }and $%
Y_{1}\rightarrow Y_{2}$ one of its elements (which must be a finitely
presentable object). For each possible strict omega functor between strict
omega categories $Y_{2}\rightarrow A$ we consider the composite morphism of
omega magmas $Y_{1}\rightarrow A$. This will be a finitely presentable
object in \textit{Pen\_Mor(A)} if $Y_{1}\rightarrow Y_{2}$ is finitely
presentable in \textit{Pen\_Mor.} Thus we have a set of finitely presentable
objects of \textit{Pen\_Mor(A)} whose elements are indexed by strict omega
functors $Y_{2}\rightarrow A$. There is one such set for each element of the
generator \textsf{G }of \textit{Pen\_Mor}. Taking the union of all these
sets over the elements of \textsf{G} gives us a set of objects that we hope
will be a strong generator for \textit{Pen\_Mor(A).} That it is indeed a
strong generator follows immediately from the fact that \textsf{G} is a
strong generator of \textit{Pen\_Mor}.

Thus we can conclude that the forgetful functor 
\[
U_{GMA}^{PMA}:\text{\textit{Pen\_Mor(A)}}\rightarrow \text{\textit{%
Omega\_Graph}}^{\rightarrow UA} 
\]

\noindent has a left adjoint $L_{PMA}^{GMA}$. Now let \textit{Omega\_Cat(A)}
denote the full subcategory of \textit{Omega\_Cat} in which an object $%
\mathbf{X}$ has the property that $X_{3}=A$. \ We define a functor 
\[
P_{A}^{\rightarrow }:\text{\textit{Omega\_Graph}}^{\rightarrow
UA}\rightarrow \text{\textit{Omega\_Cat(A)}} 
\]

\noindent as follows. Let $f:Y\rightarrow UA$ denote an object of \textit{%
Omega\_Graph}$^{\rightarrow UA}$ and denote by $\mathbf{X}$ the object $%
P_{A}^{\rightarrow }f$.

\begin{enumerate}
\item  Set $X_{1}=X_{2}\equiv domain$ $L_{PMA}^{GMA}(f)\;,\lambda ^{X}=\rho
^{X}\equiv 1_{domainL_{PMA}^{GMA}(f)}$

\item  Set $\kappa ^{X}\equiv L_{PMA}^{GMA}(f)$
\end{enumerate}

This construction is clearly functorial in $f$ and defines the functor $%
P_{A}^{\rightarrow }$.

\subsection{Categorical equivalence relations}

In applications of this theory (see part II of this work which will be
forthcoming) it is often convenient to have a criterion which identifies a
given omega magma $M$ as the domain of a (yet-to-be-constructed) categorical
Penon morphism $f$. \ Any such omega magma then determines a weak omega
category by setting $X_{1}=X_{2}=M$ \ , $\lambda ^{X}=\rho ^{X}=1_{M}$ and $%
\kappa ^{X}=f$. \ 

\begin{definition}
An \textbf{omega magma equivalence relation} on $M$ is an omega submagma $%
E\subseteq M\times M$ such that $E_{j}$ is an equivalence relation on the $j$%
-cells $M_{j}$ of $M$ for all $j\geq 0$.
\end{definition}

Recall that $U_{G}^{M}M$ is the underlying omega graph of $M$ while $%
L_{M}^{G}\circ U_{G}^{M}(M)$ is the omega magma freely generated by this
omega graph and $L_{C}^{G}\circ U_{G}^{M}(M)$ is the strict omega category
freely generated by this graph. Let \ 
\[
\varepsilon _{M}:L_{M}^{G}\circ U_{G}^{M}(M)\rightarrow M 
\]

\noindent denote the counit of the adjunction and 
\[
\eta _{M}:L_{M}^{G}\circ U_{G}^{M}(M)\rightarrow L_{C}^{G}\circ U_{G}^{M}(M) 
\]

\noindent denote the unit of the $L_{C}^{M}\dashv U_{M}^{C}$ adjunction.

\begin{definition}
An omega magma equivalence relation $E$ on $M$ is \textbf{categorical} if
for all elements $u,v\in L_{M}^{G}\circ U_{G}^{M}(M)$ the following
implication holds: $\eta _{M}u=\eta _{M}v\Longrightarrow (\varepsilon
_{M}u,\varepsilon _{M}v)\in E$.
\end{definition}

In less formal terms, we say that an omega magma equivalence relation is
categorical if any pair of elements which are equalized by \textit{every}
morphism to a strict omega category are $E$ equivalent in $M$.

\begin{definition}
An omega magma equivalence relation $E$ on $M$ is \textbf{sharp} if $a,b\in
M $ and $(cod_{k}^{j}a,dom_{k}^{j}b)\in E$ $\Longrightarrow \;\exists
\;a^{\prime },b^{\prime }\in M$ and $(a,a^{\prime }),(b,b^{\prime })\in E$
and $cod_{k}^{j}a^{\prime }=dom_{k}^{j}b^{\prime }$.
\end{definition}

\noindent Now if $E$ is any omega magma equivalence relation on $M$ we can
form the coequalizer $M//E$ of the two projections. \ This coequalizer
exists because the category of omega magmas is cocomplete. If $E$ is sharp
then this coequalizer will not contain ``extra arrows'' which in general
will arise when taking the quotient by an arbitrary equivalence relation.

\begin{theorem}
Let $E$ be a sharp, categorical, omega magma equivalence relation on $M$.
Then $(M//E)_{j}\backsimeq M_{j}/E_{j}$ and $M//E$ is a strict omega
category.
\end{theorem}

\noindent \textbf{Proof:}

It will suffice to show that the coequalizer of the two projections from $%
U_{G}^{M}E$ in the category of omega graphs can be given the structure of a
strict omega category in such a way that the quotient map from $U_{G}^{M}M$
to the coequalizer is in fact an omega magma morphism. If this can be done
then this strict omega category is easily seen to have the universal
property of the desired coequalizer and hence must be isomorphic to it.

Now the coequalizer $U_{G}^{M}E\rightrightarrows U_{G}^{M}M$ is just the
omega graph

\noindent $U_{G}^{M}M//U_{E}^{M}E$ whose $j$-cells are the equivalence
classes of $M_{j}$ by the equivalence relation $E_{j}$; this set of
equivalence classes was denoted in the statement of the theorem by $%
M_{j}/E_{j}$. We now show that $U_{G}^{M}M//U_{E}^{M}E$ has an omega magma
structure for which the quotient map is an omega magma morphism.

Let $[a],[b]$ be equivalence classes in $M_{j}/E_{j}$ such that $%
cod_{k}^{j}[a]=dom_{k}^{j}[b]$. Since $U_{G}^{M}M//U_{E}^{M}E$ is an omega
graph we see that $(cod_{k}^{j}a,dom_{k}^{j}b)\in E$. Since $E$ \ is sharp
there exists $a^{\prime }\in \lbrack a],b^{\prime }\in \lbrack b]$ such that 
$cod_{k}^{j}a^{\prime }=dom_{k}^{j}b^{\prime }$. Then define $[a]\bigcirc
_{k}^{j}[b]\equiv \lbrack a^{\prime }\bigcirc _{k}^{j}b^{\prime }]$. Because 
$E$ is an omega submagma of $M$ this definition is independent of the
choices made. Clearly the quotient map $U_{G}^{M}M\rightarrow
U_{G}^{M}M//U_{E}^{M}E$ \ is the underlying omega graph morphism of a
morphism of omega magmas.

The fact that the omega magma structure defined on $U_{G}^{M}M//U_{E}^{M}E$
gives it the structure of a strict omega category follows immediately from
the fact that the relation $E$ was assumed to be categorical.$\blacksquare $

Now let us assume that the omega magma $M$ comes equipped with a bridge
relation $R^{M}$ and a sharp, categorical, omega magma equivalence relation $%
E$.

\begin{definition}
$R^{M}$ is a \textbf{witness} to $E$ if
\end{definition}

\begin{enumerate}
\item  $(a,b,c)\in R_{j}^{M}\Longrightarrow (a,b)\in E$ and $\
(c,id_{j+1}^{j}a),(c,id_{j+1}^{j}b)\in E$

\item  $a\parallel b$ and $(a,b)\in E_{j}\Longrightarrow \exists c$ and $%
(a,b,c)\in R_{j}^{M}$
\end{enumerate}

\noindent The following result is a trivial consequence of the last
definition and of the preceding theorem.

\begin{theorem}
Let $(M,R^{M})$ be a bridge magma, let $E$ be a sharp, categorical, omega
magma equivalence relation on $M$ and assume $R^{M}$ is a witness to $E$.
Then the quotient map induces a categorical Penon morphism. 
\[
\kappa _{E}:(M,R^{M})\rightarrow (M//E,R^{\Delta }) 
\]
\end{theorem}

\section{Pseudo-functors}

In this section we develop a theory of omega pseudo-functors between weak
omega categories. Intuitively, an omega pseudo-functor differs from an omega
functor in that the former preserves operations and identities only up to
equivalence instead of ``on the nose''. It is our view that omega
pseudo-functors are the most natural type of morphism between weak omega
categories. We show that omega pseudo-functors are ubiquitous by showing
that in any weak omega category $\mathbf{X}$ composition with an $i$-cell
defines a pair of omega pseudo-functors that are not in general omega
functors. More specifically, we show that if $a,b,c$ is a triple of $\left(
i-1\right) $-cells in $X_{1}$ and if $h$ is an $i$-cell with $%
dom_{i-1}^{i}h=b$ and $cod_{i-1}^{i}h=c$ then there is an omega
pseudo-functor 
\[
\mathbf{\Theta (h):X(}a,b)\rightarrow \mathbf{X(}a,c) 
\]

\noindent which is defined by right composition with $h$ over the $\left(
i-1\right) $-cell $b$ . Of course the analogous result holds for left
composition as well.

We also show that an omega pseudo-functor between 2-skeletal weak omega
categories is a homomorphism in the standard sense (which we shall call
``classical''). The classical homomorphisms which arise in this way from
omega pseudo-functors have a special property of being ``proper''. We know
of no classical homomorphism that is not also proper but we expect that such
examples can be constructed. In any case non-proper homomorphisms do not
seem to arise naturally. We close this section by showing that any proper,
classical pseudo-functor is an omega pseudo-functor.

The reader will notice that to establish the correspondence between omega
pseudo-functors and classical, proper homomorphisms in the 2-skeletal case
we have to employ some lengthy arguments. This suggests to us that our
definition of pseudo-functor is indeed a non-trivial change of viewpoint
even for homomorphisms between bicategories.

\subsection{Definition of omega pseudo-functor}

Let $\mathbf{X,Y}$ be weak omega categories and let $\mathbf{F}%
=(F_{1},F_{2},F_{3})$ be a triple of morphisms of \textit{globular sets.}%
Consider the following diagram of globular sets: 
\[
\begin{array}{cccccc}
&  & {\small \lambda }^{X} &  & {\small \kappa }^{X} &  \\ 
& X_{1} & \leftarrow & X_{2} & \rightarrow & X_{3} \\ 
& \mid &  & \mid &  & \mid \\ 
{\small F}_{1} & \mid & {\small F}_{2} & \mid & {\small F}_{3} & \mid \\ 
& \mid &  & \mid &  & \mid \\ 
& \downarrow & {\small \lambda }^{Y} & \downarrow & \kappa ^{Y} & \downarrow
\\ 
& Y_{1} & \longleftarrow & Y_{2} & \rightarrow & Y_{3}
\end{array}
\]
\bigskip\ 

\begin{definition}
$\mathbf{F}$ is an \textbf{omega pseudo-functor }if

\begin{enumerate}
\item  $F_{3}$ is a morphism of omega graphs

\item  the $F_{i}$ make the above diagram commute in the category of
globular sets

\item  for all $a,b$ such that $a\bigcirc _{j-1}^{j}b$ is defined in $X_{2}$ 
\[
\kappa ^{Y}(F_{2}a\bigcirc _{j-1}^{j}F_{2}b)=\kappa ^{Y}F_{2}(a\bigcirc
_{j-1}^{j}b) 
\]
\end{enumerate}
\end{definition}

We note that we do \textbf{not} require the equality $F_{2}\circ \rho
^{X}=\rho ^{Y}\circ F_{1}$ to hold for $\mathbf{F}$ to be an omega
pseudo-functor.

\begin{proposition}
Let $\mathbf{F:X}\rightarrow \mathbf{Y}$ and $\mathbf{G:Y}\rightarrow 
\mathbf{Z}$ be omega pseudo-functors. Then $\mathbf{G\circ F}$ is an omega
pseudo-functor.
\end{proposition}

\noindent \textbf{Proof:}

The composite $G_{3}\circ F_{3}$ is a morphism of omega graphs so it remains
only to verify condition 3 of the definition. To this end let $a,b\in X_{2}$
and suppose that $a\bigcirc _{j-1}^{j}b$ is defined. Since $\mathbf{F}$ is
an omega pseudo-functor we know that 
\[
\kappa ^{Y}(F_{2}a\bigcirc _{j-1}^{j}F_{2}b)=\kappa ^{Y}\circ
F_{2}(a\bigcirc _{j-1}^{j}b) 
\]

\noindent Since $G_{3}\circ \kappa ^{Y}=\kappa ^{Z}\circ G_{2}$ we conclude 
\[
\kappa ^{Z}\circ G_{2}(F_{2}a\bigcirc _{j-1}^{j}F_{2}b)=\kappa ^{Z}\circ
G_{2}\circ F_{2}(a\bigcirc _{j-1}^{j}b)\text{.} 
\]

\noindent On the other hand the fact that $\mathbf{G}$ is an omega pseudo
functor yields 
\[
\kappa ^{Z}((G_{2}\circ F_{2}a)\bigcirc _{j-1}^{j}(G_{2}\circ
F_{2}b))=\kappa ^{Z}\circ G_{2}(F_{2}a\bigcirc _{j-1}^{j}F_{2}b)\text{.} 
\]

\noindent Therefore 
\[
\kappa ^{Z}((G_{2}\circ F_{2}a)\bigcirc _{j-1}^{j}(G_{2}\circ
F_{2}b))=\kappa ^{Z}\circ G_{2}\circ F_{2}(a\bigcirc _{j-1}^{j}b) 
\]

\noindent as desired. $\blacksquare $

From the preceding proposition we conclude that there is a category \textit{%
PF\_Omega\_Cat} in which an object is a weak omega category and a morphism
is an omega pseudo-functor. \textit{PF\_Omega\_Cat }contains \textit{%
Omega\_Cat} as a subcategory. As a category \textit{PF\_Omega\_Cat }is not
nearly so well-behaved as \textit{Omega\_Cat}. For example, it is neither
complete nor cocomplete. \ Thus it is not locally presentable and is not the
category of models of an essentially algebraic theory.

We note one more obvious fact. Let $\mathbf{F:X\rightarrow Y}$ be an omega
pseudo-functor. Let $a,b$ be any pair of $(i-1)$-cells of $X_{1}$. Then the
appropriate restriction of $\mathbf{F}$ defines an omega pseudo-functor 
\[
\mathbf{F(}a,b)\mathbf{:X}(a,b)\mathbf{\rightarrow Y}(F_{1}a,F_{1}b) 
\]

\subsection{\protect\bigskip Pseudo-functors defined by composition in a
weak omega category}

We next show that every $i$-cell $h$, $i\geq 1$, in the underlying omega
magma of a weak omega category $\mathbf{X}$ defines a pair of omega
pseudo-functors which correspond to left and right composition with $h$. We
shall prove this for right composition since the case of left composition is
proved in exactly the same way.

Let $a,b,c$ be three $\left( i-1\right) $-cells of $X_{1}$and $h$ an $i$%
-cell satisfying 
\[
dom_{i-1}^{i}h=b,\;cod_{i-1}^{i}h=c 
\]
Recall the construction of the weak omega categories $\mathbf{X}(a,b),%
\mathbf{X}(a,c)$ from Section 7. We shall define a triple of globular set
morphisms $\Theta (h)_{j},\;j=1,2,3$ with 
\[
\Theta (h)_{j}:X(a,b)_{j}\rightarrow X(a,c)_{j} 
\]

\noindent and then we shall show that $\mathbf{\Theta (h)\equiv (}\Theta
(h)_{j})$ is an omega pseudo-functor. Let $\bigcirc $ denote composition in
any of the three omega magmas defining $\mathbf{X}$.

\begin{enumerate}
\item  Let $f\in (X(a,b)_{1})_{j}$. By definition this means that $f\in
(X_{1})_{i+j}$. Define 
\[
\Theta (h)_{1}f\equiv f\bigcirc _{i-1}^{i+j}id_{i+j}^{i}h 
\]

\item  Let $f\in (X(a,b)_{2})_{j}$. Then $f\in (X_{2})_{i+j}$. Define 
\[
\Theta (h)_{2}f\equiv f\bigcirc _{i-1}^{i+j}id_{i+j}^{i}(\rho ^{X}h) 
\]

\item  Let $f\in (X(a,b)_{3})_{j}$. Then $f\in (X_{3})_{i+j}$. Define 
\[
\Theta (h)_{3}f\equiv f\bigcirc _{i-1}^{i+j}id_{i+j}^{i}(\kappa ^{X}\rho
^{X}h) 
\]
\end{enumerate}

We note that since $X_{3}$ is a strict omega category the morphism of
globular sets $\Theta (h)_{3}$ is in fact a strict functor between strict
omega categories because the interchange law and identity laws hold as
equations in $X_{3}$. We shall see that this is a special property of omega
pseudo-functors of the form $\mathbf{\Theta (h)}$ which is not shared by the
general omega pseudo-functor.

\begin{proposition}
$\mathbf{\Theta (h)}:X(a,b)\rightarrow X(a,c)$ is an omega pseudo-functor.
\end{proposition}

\noindent \textbf{Proof:}

One easily sees that the triple of globular set morphisms

\noindent $(\Theta (h)_{1},\Theta (h)_{2},\Theta (h)_{3})$ makes the
relevant diagram of globular sets commute. Since $\Theta (h)_{3}$ is a
morphism of strict omega categories it is a fortiori a morphism of omega
graphs. It thus only remains to check condition 3 in the definition of omega
pseudo-functor.

Let $f,g\in (X(a,b)_{2})_{j}$ and suppose that $f\bigcirc _{j-1}^{j}g$ is
defined. Then $f,g\in (X_{2})_{i+j}$ and $f\bigcirc _{i+j-1}^{i+j}g$ is
defined and 
\[
\Theta (h)_{2}(f\bigcirc _{j-1}^{j}g)\equiv (f\bigcirc
_{i+j-1}^{i+j}g)\bigcirc _{i-1}^{i+j}id_{i+j}^{i}h 
\]

\noindent where the left hand side of this equation is an expression
involving operations and elements of $X(a,b)_{2}$ while the right hand side
is an expression involving operations and elements of $X_{2}$. Define 
\[
u\equiv (f\bigcirc _{i+j-1}^{i+j}g)\bigcirc _{i-1}^{i+j}id_{i+j}^{i}h\text{.}
\]

\noindent Note that

\[
\Theta (h)_{2}(f)\bigcirc _{j-1}^{j}\Theta (h)_{2}(g)=(f\bigcirc
_{i-1}^{i+j}id_{i+j}^{i}h)\bigcirc _{i+j-1}^{i+j}(g\bigcirc
_{i-1}^{i+j}id_{i+j}^{i}h) 
\]

\noindent where again the left hand side denotes operations and elements in

\noindent $X(a,b)_{2}$ while the right hand side denotes operations and
elements in $X_{2}$. Define 
\[
v\equiv (f\bigcirc _{i-1}^{i+j}id_{i+j}^{i}h)\bigcirc
_{i+j-1}^{i+j}(g\bigcirc _{i-1}^{i+j}id_{i+j}^{i}h)\text{.} 
\]

\noindent We must show that $\kappa ^{X}u=\kappa ^{X}v$.\ 

To this end we define an element $u^{\prime }\,$\ of $(X_{2})_{i+j}$ 
\[
u^{\prime }\equiv (f\bigcirc _{i+j-1}^{i+j}g)\bigcirc
_{i-1}^{i+j}(id_{i+j}^{i}h\bigcirc _{i+j-1}^{i+j}id_{i+j}^{i}h)\text{.} 
\]

\noindent Note that since identity laws hold strictly in $X_{3}$ and since $%
\kappa ^{X}$ is an omega magma morphism $\kappa ^{X}u=\kappa ^{X}u^{\prime }$%
. But the interchange law also holds as an equation in $X_{3}$ so $\kappa
^{X}u^{\prime }=\kappa ^{X}v$. $\blacksquare $

\subsection{\protect\bigskip Omega pseudo-functors and homomorphisms}

We shall now examine the relationship between omega pseudo-functors and
classical pseudo-functors (usually called homomorphism in the bicategory
case and pseudo-functors in the strict 2-category case).\ Recall that the
results of Section 4 show that the category of bicategories and strong
homomorphisms is a retract of the category of 2-skeletal, weak omega
categories and omega functors. We cannot prove an analogous result in which
homomorphisms replace strong homomorphism and omega-pseudo functors replace
omega functors. In fact we strongly suspect \ that such a retraction may not
exist.

In any event what we will show is that the $F_{1}$ component of any omega
pseudo functor between 2-skeletal, weak omega categories $\mathbf{X}$ and $%
\mathbf{Y}$ defines a classical homorphism between the bicategories $X_{1}$
and $Y_{1}$. The bicategory homomorphisms which arise in this way have a
special property we call ``properness''. We have been unable to construct or
to find any example of a homomorphism that is not proper. However, the
definition of homomorphism does not apparently preclude the existence of
such an example.

In the next subsection we shall also show that for any classical
homomorphism $h$ between bicategories $A$ and $B$ there exist 2-skeletal,
weak omega categories $\mathbf{X}$ and $\mathbf{Y}$ and omega pseudo-functor 
$\mathbf{F}$ between them such that $X_{1}=A$, $Y_{1}=B$ and $F_{1}=h$. The
proof of this last result is surprisingly lengthy and is evidence that our
definition of omega pseudo-functor represents a significant change of
viewpoint from the classical one.

Let $A,B$ be classical bicategories. We refer to [7,12,14] for the
definition of a homorphism from $A$ to $B$. Essentially a homomorphism
differs from a strong homomorphism in that it preserves identities and
operations only up to isomorphisms in $B$. These isomorphisms must fit
together ``coherently'' and must thus satisfy some axioms. These axioms make
working with composite homomorphisms a complicated task. Our definition of
omega pseudo-functor offers a path through these complications.

Let $H:A\rightarrow B$ be a homorphism between classical bicategories. Part
of the definition of $H$ is coherence data which must satisfy some axioms.
This coherence data has two components:

\begin{enumerate}
\item  For each pair of 1-cells $f,g\in A$ such that $f\bigcirc _{0}^{1}g$
is defined there is given a 2-cell isomorphism 
\[
\phi _{f,g}:Hf\bigcirc _{0}^{1}Hg\rightarrow H(f\bigcirc _{0}^{1}g) 
\]

\item  For each 0-cell $a\in A$ there is given a 2-cell isomorphism 
\[
\phi _{a}:id_{1}^{0}(Ha)\rightarrow H(id_{1}^{0}a). 
\]
\end{enumerate}

\noindent The 2-cell isomorphisms $\phi _{f,g}$ and $\phi _{a}$ are in fact
components of a pair of natural isomorphisms between functors. These natural
transformations must satisfy some axioms which we shall not display here
(but see [7,12].

\begin{definition}
\noindent A homomorphism $H:A\rightarrow B$ between classical bicategories
is called \textbf{proper} if the associated coherence data has the following
properties:

\begin{enumerate}
\item  for all 1-cells $f,g$ of $A$ 
\[
Hf\bigcirc _{0}^{1}Hg=H(f\bigcirc _{0}^{1}g)\Longrightarrow \phi
_{f,g}=id_{2}^{1}(H(f\bigcirc _{0}^{1}g)) 
\]

\item  for all 0-cells $a$ of $A$ 
\[
id_{1}^{0}(Ha)=H(id_{1}^{0}a)\Longrightarrow \phi
_{a}=id_{2}^{1}(H(id_{1}^{0}a)) 
\]

\item  for all 1-cells $f,f^{\prime },g,g^{\prime }$ of $A$ 
\[
f\bigcirc _{0}^{1}g=f^{\prime }\bigcirc _{0}^{1}g^{\prime }\text{ and \ }%
Hf\bigcirc _{0}^{1}Hg=Hf^{\prime }\bigcirc _{0}^{1}Hg^{\prime
}\Longrightarrow \phi _{f,g}=\phi _{f^{\prime },g^{\prime }} 
\]

\item  for all 0-cells $a,b$ of $A$ $\ $%
\[
Ha=Hb\text{ and }H(id_{1}^{0}a)=H(id_{1}^{0}a)\Longrightarrow \phi _{a}=\phi
_{b} 
\]
\end{enumerate}
\end{definition}

\begin{theorem}
Let $\mathbf{F:X\rightarrow Y}$ be an omega pseudo-functor between
2-skeletal, weak omega categories. Then there is coherence data $\phi _{f,g}$
and $\phi _{a}$ which exhibits $F_{1}$ as a homomorphism between the
bicategories $X_{1}$ and $Y_{1}$. This homomorphism is proper.
\end{theorem}

\noindent \textbf{Proof:}

The fact that $X_{1}$ and $Y_{1}$ are bicategories was proven in Section 4.
We shall construct the coherence data $\phi _{f,g}$ and $\phi _{a}$ and then
outline the method the reader may follow to satisfy himself that all
necessary diagrams involving this coherence data commute in $Y_{1}$. It will
be clear from the construction of this coherence data that $F_{1}$ must be a
proper homomorphism.

Let $f,g$ be 1-cells of $X_{1}$ such that $f\bigcirc _{0}^{1}g$ is defined.
Since $\mathbf{F}$ is an omega pseudo-functor we know that $\kappa ^{Y}$
equalizes the 1-cells $u,v$ of $Y_{2}$ where these cells are defined by the
equations 
\[
u\equiv (F_{2}\circ \rho ^{X})f\bigcirc _{0}^{1}(F_{2}\circ \rho ^{X})g\; 
\]
\[
v\equiv (F_{2}\circ \rho ^{X})(f\bigcirc _{0}^{1}g). 
\]

\noindent Since these two cells are also parallel we conclude that there
exists a (unique) 2-cell $\phi _{f,g}^{\prime }\in Y_{2}$ such that 
\[
(u,v,\phi _{f,g}^{\prime })\in R_{1}^{Y}. 
\]
(Later in the proof it will be helpful for the reader to recall that $\kappa
^{Y}$ $\phi _{f,g}^{\prime }$ is an identity 2-cell of $Y_{3}$.) \ Define $%
\phi _{f,g}\equiv \lambda ^{Y}\phi _{f,g}^{\prime }$. \ Since $\lambda ^{Y}$
is an omega magma morphism that is split by the omega graph morphism $\rho
^{Y}$ and since $\lambda ^{Y}\circ F_{2}=F_{1}\circ \lambda ^{X}$ as
morphisms of globular sets we conclude that 
\[
dom_{1}^{2}\phi _{f,g}=F_{1}f\bigcirc _{0}^{1}F_{1}g 
\]

\noindent and 
\[
cod_{1}^{2}\phi _{f,g}=F_{1}(f\bigcirc _{0}^{1}g) 
\]

We define the coherence data $\phi _{a}$ in a similar way. Here the
assumption that $F_{3}$ is a morphism not just of globular sets but of omega
graphs (and hence preserves identities) will play a crucial role.

Let $a$ denote a 0-cell of $X_{1}$ and consider the 1-cells of $Y_{2}$
defined by the equations 
\[
u\equiv F_{2}\,id_{1}^{0}(\rho ^{X}a) 
\]
\[
v\equiv id_{1}^{0}(F_{2}\rho ^{X}a). 
\]

\noindent These two 1-cells are obviously parallel. We claim that they are
equalized by $\kappa ^{Y}$. \ To see this first note that by hypothesis $%
\kappa ^{Y}\circ F_{2}=F_{3}\circ \kappa ^{X}$. \ The latter morphism is
actually a morphism of omega graphs because $F_{3}$ is assumed to be a
morphism of omega graphs. Thus $\kappa ^{Y}\circ F_{2}$ is a morphism of
omega graphs. It follows that 
\[
\kappa ^{Y}u=\kappa ^{Y}F_{2}(id_{1}^{0}(\rho ^{X}a))=id_{1}^{0}(\kappa
^{Y}F_{2}\rho ^{X}a)=\kappa ^{Y}id_{1}^{0}(F_{2}\rho ^{X}a)=\kappa ^{Y}v. 
\]

\noindent Consequently we deduce the existence of a (unique) 2-cell $\phi
_{a}^{\prime }$ of $Y_{2}$ such that 
\[
(u,v,\phi _{a}^{\prime })\in R_{1}^{Y}. 
\]

\noindent We define $\phi _{a}\equiv \lambda ^{Y}\phi _{a}^{\prime }$ and
observe as before that 
\[
dom_{1}^{2}\phi _{a}=F_{1}id_{1}^{0}(a) 
\]
\[
cod_{1}^{2}\phi _{a}=id_{1}^{0}(F_{1}a). 
\]

It remains to perform the straightforward but laborious task of verifying
that these coherence data define the required natural transformations and
that these natural transformations satisfy the required axioms. This in turn
amounts to showing that certain diagrams of 2-cells in $Y_{1}$ commute (see
[7,12] for details.) \ These verifications all follow the same pattern which
was outlined in the Scholium of Section 3. For the reader's convenience we
shall again describe this method here.

The typical diagram of 2-cells in $Y_{1}$ which we must show is commutative
involves two distinct composite 2-cells which a priori are parallel in $%
Y_{1} $. \ To show they are in fact equal it suffices to exhibit a 3-cell of 
$Y_{1} $ which is a bridge from one of these composite cells to the other.
For by hypothesis all cells of $Y_{1}$ above dimension 2 are identity cells.

To construct the desired bridge cells one first ``disassembles'' the diagram
in $Y_{1}$. Note that it consists of images under $F_{1}$ of 1-cells of $%
X_{1}$ together with coherence 2-cells for $F_{1}$, and possibly the image
under $F_{1}$ of the associator of $X_{1}$ and the associator for $Y_{2}$. \
One ``reassembles'' this diagram in $Y_{2}$ by first sending the individual
1-cells of $X_{1}$ to $X_{2}$ via $\rho ^{X}$, composing them as appropriate
and then sending them to $Y_{2}$ via $F_{2}$. These cells are then
reconnected by the 2-cells of $Y_{2}$ whose images under $\lambda ^{Y}$ and $%
\lambda ^{X}$ define the coherence cells for $F_{1}$ and the coherence data
for $X_{1}$ and $Y_{1}$.

Once the diagram has be reassembled in $Y_{2}$ one simply notes that the
images under $\kappa ^{Y}$ of each of its 2-cells is an identity in $Y_{3}$
because of the definition of the coherence cells as images of bridge cells
in $Y_{2}$. It follows that the two composite 2-cells in the diagram in $%
Y_{2}$ are bridged in $Y_{2}$ and the image under $\lambda ^{Y}$ of this
bridge cell is the desired bridge cell in $X_{1}$. $\blacksquare $

\subsection{\protect\bigskip Homomorphisms are omega pseudo-functors}

We next turn our attention to the following question. Suppose 
\[
h:A\rightarrow B 
\]
is a homomorphism between classical bicategories. Is there an omega
pseudo-functor between 2-skeletal, weak omega categories whose first
coordinate is $h$?

\begin{theorem}
Let $h:A\rightarrow B$ be a homomorphism between classical bicategories. If $%
h$ is proper then there exists an omega pseudo-functor $\mathbf{%
F:X\rightarrow Y}$ where $\mathbf{X,Y}$ are 2-skeletal, $X_{1}=A$, $Y_{1}=B$
and $F_{1}=h$.
\end{theorem}

The proof of this result is an unpleasantly long construction but it is
entirely straight forward. We shall outline the steps so that the reader may
satisfy herself that the details can be filled in as necessary.

In Section 4 we constructed a functor 
\[
S_{C}^{B}:\text{\textit{Bicat}}\rightarrow \text{\textit{Weak\_2Cat.}} 
\]

\noindent This functor takes a bicategory $A$ to a 2-skeletal, weak omega
category $\mathbf{X}$ whose $X_{1}$ component is $A$, whose $X_{2}$
component is the bicategory freely generated by the omega graph underlying $%
A $, and whose $X_{3}$ component is the strict 2-category freely generated
by the omega graph underlying $A$.

\begin{enumerate}
\item  Define $\mathbf{X}=$\ $S_{C}^{B}A$ and $\mathbf{Z}=S_{C}^{B}B$. \
Define $F_{1}:X_{1}\rightarrow Y_{1}$ to be $h$. Define $F_{2}:X_{2}%
\rightarrow Y_{2}$ by the equation 
\[
F_{2}\equiv \rho ^{Y}\circ F_{1}\circ \lambda ^{X} 
\]
The definition of $F_{3}$ will be left to the last step of the construction.
(We cannot simply use $F_{1}$ to define $F_{3}$ because the result will not
be a map of omega graphs.) \ The plan is to modify $Z_{2}$ by adding
appropriate 2-cell equivalences and 3-cell isomorphisms and to modify $Z_{3}$
by taking a quotient by the appropriate equivalence relation so that the
result is a 2-skeletal, weak omega category $\mathbf{Y}$ with the desired
properties. These properties will allow the construction of $F_{3}$ which
will be an omega graph morphism to $Y_{3}$ that does not in general factor
through $Z_{3}$.

\item  Define an omega magma equivalence relation (definition 24) $E$ on $%
Z_{2}$ to be the smallest such relation generated by all pairs of the form:

\begin{enumerate}
\item  $(a,b)$ where $(a,b,c)\in R^{Z}$

\item  $(F_{2}(f\bigcirc _{0}^{1}g),F_{2}f\bigcirc _{0}^{1}F_{2}g)$

\item  $(F_{2}(id_{1}^{0}a),id_{1}^{0}(F_{2}a))$
\end{enumerate}

For later purposes it is important to note that $(a,b)\in
E_{1}\Longrightarrow a\parallel b$. This follows from the fact that $E_{1}$
is generated by parallel pairs. Moreover, if $(a,b)\in E$ and $a\parallel b$
then either $a=b$ or $(a,b)\in E_{1}$.

\item  Construct a function $\psi :E_{1}\rightarrow (Z_{1})_{2}$ which takes
a pair to its associated coherence 2-cell isomorphism in $Z_{1}$. The main
difficulty here is to guarantee that $\psi $ is well defined. The
construction of $\psi $ proceeds as follows. Let $D_{1}$ denote the
equivalence relation induced on $Z_{2}$ by the ternary relation $R_{1}^{Z}$.
One first observes that $D_{1}$ is closed in $E_{1}$ under the operation $%
\bigcirc _{0}^{1}$. This last assertion follows from the coherence theorem
for bicategories. Clearly $D_{1}$ also contains the diagonal of $E_{1}$. We
now observe that $E_{1}-D_{1}$ is a 2-sided ideal in $E_{1}$ in the sense
that left or right composition with an element of $D_{1}$ again yields an
element of $E_{1}-D_{1}$. This observation is justified as follows. Since $%
Z_{3}$ is freely generated by the same omega graph which underlies $Z_{2}$
we know that 
\[
(a.b)\in E_{1}-D_{1}\Leftrightarrow \kappa ^{Z}a\neq \kappa ^{Z}b\text{.} 
\]

Since pairs in $D_{1}$ are equalized by $\kappa ^{Z}$ and because $Z_{3}$ is
freely generated we conclude that 
\[
(a,b)\in E_{1}-D_{1}\text{ and }(d,d^{\prime })\in D_{1}\Longrightarrow
\kappa ^{Z}(a\bigcirc _{0}^{1}d)\neq \kappa ^{Z}(b\bigcirc _{0}^{1}d^{\prime
}) 
\]

and similarly for left composition.

The fact that $E_{1}-D_{1}$ is a 2-sided ideal in $E_{1}$ allows us to
separate the construction of $\psi $ into two parts. The function 
\[
\psi \mid D_{1}\equiv \lambda ^{Z}\mid D_{1}. 
\]
The function $\psi \mid E_{1}-D_{1}$ is defined by induction. Recall that $%
Z_{2}$ is freely generated by the omega graph underlying the bicategory $B$.
\ We associate with a pair of $E_{1}-D_{1}$ the integer which is the minimum
of the word lengths of the two elements of the pair. We begin the inductive
construction for $n=1$ by defining $\psi $ using the coherence data for $h$.
For the inductive step, we observe that any pair $(a,b)$ of length $n$ can
be written in at least one way as the composite of two pairs of lengths $%
\leq n-1$ one of which must be in $E_{1}-D_{1}$: 
\[
(a,b)=(a^{\prime },b^{\prime })\bigcirc _{0}^{1}(a^{\prime \prime
},b^{\prime \prime }) 
\]
Now we observe that the fact that the homomorphism $h$ is proper and the
fact that its classical coherence data for $h$ satisfies the classical
axioms means that the assignment: 
\[
\psi (a,b)\equiv \psi (a^{\prime },b^{\prime })\bigcirc _{0}^{2}\psi
(a^{\prime \prime },b^{\prime \prime }) 
\]

is independent of the choice of the decomposition of $(a,b)$.

\item  Observe that since $\kappa ^{Z}$ is a surjection in the category of
omega graphs $\kappa ^{Z}E$ is an omega magma equivalence relation on $Z_{3}$%
. Define $Y_{3}$ to be the coequalizer of this relation in the category of
strict omega categories and let $\widehat{\kappa }^{Z}$ denote the composite
of $\kappa ^{Z}$ followed by the map to the coequalizer $Y_{3}$.

\item  Observe that 
\[
a,b\in Z_{2}\text{ and }\widehat{\kappa }^{Z}a=\widehat{\kappa }^{Z}b\text{
and }a\parallel b\Longrightarrow (a,b)\in E. 
\]
This is true for pairs already equalized by $\kappa ^{Z}$ by definition of $%
E $ and for other pairs by definition of the coequalizer and the fact that $%
\kappa ^{Z}E$ is transitive.

\item  Construct $Y_{2}$, a 3-skeletal omega magma, in two stages.

In the first stage define $Z_{2}^{\prime }$ to be the bicategory obtained
from $Z_{2}$ by freely adjoining one 2-cell for each element $(a,b)$ of $%
E_{1}$. (Recall that $Z_{2}$ is itself a free bicategory.) The domain of
this 2-cell is $a$ an its codomain is $b$. This gives us a monic morphism of
bicategories 
\[
\chi :Z_{2}\rightarrow Z_{2}^{\prime } 
\]
which is an isomorphism on 0- and 1-cells. Define $\lambda ^{Z^{\prime }}$
as the extension of $\lambda ^{Z}$ obtained by sending the formal 2-cell
associated with the element $(a,b)$ to the 2-cell $\psi (a,b)$ of $Z_{1}$
defined in step 3 and extending multiplicatively. Extend $\widehat{\kappa }%
^{Z}$ over these new 2-cells by sending them to the appropriate identity
cells of $Y_{3}$.

In the second stage add a unique 3-cell isomorphism to $Z_{2}^{\prime }$ for
each parallel pair of 2-cells of $Z_{2}^{\prime }$ equalized by $\widehat{%
\kappa }^{Z}$. \ This defines a 3-skeletal omega magma $Z_{2}^{\prime \prime
}$. Then define $\lambda ^{Z^{\prime \prime }}$ by extending $\lambda
^{Z^{\prime }}$ over $Z_{2}^{\prime \prime }$ by mapping these 3-cells to
identities in $Z_{3}$. Define $\widehat{\kappa }^{Z^{\prime \prime }}$ by
extending $\widehat{\kappa }^{Z}$ over these new 3-cells by mapping them to
identities in $Y_{3}$.

\item  Define $\mathbf{Y}$ by setting $Y_{1}=Z_{1}$, $Y_{2}=Z_{2}^{\prime
\prime }$, $\lambda ^{Y}=\lambda ^{Z^{\prime \prime }}$, $\rho ^{Y}=\rho
^{Z} $, $Y_{3}$ defined in step 4 and $\kappa ^{Y}=\widehat{\kappa }%
^{Z^{\prime \prime }}$.

\item  Observe that $\mathbf{Y}$ is a 2-skeletal, weak omega category
because $\kappa ^{Y}$ is a categorical Penon morphism by construction.
Observe that the map of globular sets 
\[
F_{2}:X_{2}\rightarrow Z_{2}\hookrightarrow Y_{2} 
\]
has the properties required of the second coordinate of an omega
pseudo-functor by construction of $Y_{3}$.

\item  Finally we construct a morphism of omega graphs $F_{3}:X_{3}%
\rightarrow Y_{3}$ which will fill in the necessary commutative diagram.
This is easily done by observing that if $a\in X_{3}$ we may choose any
element $b\in X_{2}$ such that $\kappa ^{X}b=a$ and then define $%
F_{3}a\equiv \kappa ^{Y}\circ F_{2}b$. This is well defined for the
following reason. Observe that because $X_{2}$ and $X_{3}$ are freely
generated $\kappa ^{X}$ equalizes only those parallel pairs of elements that
must be equalized by any omega magma morphism to any strict omega category.
For any parallel pair of elements $(b,b^{\prime })\in X_{2}$ which are
equalized by $\kappa ^{X}$ we can find a parallel pair of element $\left(
c,c^{\prime }\right) \in Y_{2}$ which are equalized by $\kappa ^{Y}$ and
such that $(F_{2}b,c)$, $(F_{2}b^{\prime },c^{\prime })\in \dot{E}$. This
follows from the last observation and from the definitions of $F_{2}$ and $E$%
. Consequently $(F_{2}b,F_{2}b^{\prime })\in E$ and are therefore equalized
by $\kappa ^{Y}$.

The morphism $F_{3}$ defined in this way is in fact a morphism of omega
graphs by the same argument.$\blacksquare $
\end{enumerate}

\section{Omega equivalence of weak omega categories}

We think it fair to assert that standard (1 dimensional) category theory is
built around the notion of isomorphism between objects and the attendant
concepts of universal arrow and adjunction. Two dimensional category theory
permits a coarsening of the isomorphism relation via the notion of a 1-cell
equivalence between objects. This leads to more general notions of limit
(weighted limits). In higher dimensions one would expect that theories of
higher dimensional limits and adjunctions would depend upon suitable notions
of higher dimensional equivalence which would be still coarser than the
notion of equivalence in the two dimensional case.

In this section we develop two good candidates for this higher dimensional
notion of equivalence between weak omega categories.\ The coarser one, 
\textit{weak equivalence}, generalizes the ordinary notion of categorical
equivalence and biequivalence between bicategories. It is akin to the notion
of weak homotopy equivalence in homotopy theory and is induced by a single
omega functor or omega pseudo-functor. The stronger relation, \textit{omega
equivalence}, differs from weak equivalence only in that it is defined by
two omega functors or omega pseudo-functors in opposite directions. \ In
this aspect it is similar to the notion of homotopy equivalence in homotopy
theory.

Both notions of equivalence depend upon the construction of a functor 
\[
\Pi :\text{\textit{Tame\_Omega\_Cat}}\rightarrow \text{\textit{Set}} 
\]

\noindent This functor is a generalization of the functor that associates
with an ordinary category its set of isomorphism classes. The category 
\textit{Tame\_Omega\_Cat} has the same objects as \textit{PF\_Omega\_Cat, }%
the category in which an object is a weak omega category and a morphism an
omega pseudo-functor. The morphisms of \textit{Tame\_Omega\_Cat} are what we
shall call \textit{tame} omega pseudo-functors. An omega pseudo-functor is
tame if it preserves the cells (called internal equivalences) that are used
to define the functor $\Pi $ on objects. For example if $\mathbf{X}$ is a
1-skeletal weak omega category then $\Pi \mathbf{X}$ is the set of
isomorphism components of $X_{1}$ (which was shown to be an ordinary
category in Section 3). If $\mathbf{X}$ is a 2-skeletal weak omega category
then $\Pi \mathbf{X}$ is the set if equivalence classes of objects of $X_{1}$
(a bicategory by Section 4) where two objects are equivalent if they are
connected by a 1-cell equivalence.

The notions of weak equivalence and omega equivalence are the cornerstones
of a theory of weighted limits for weak omega categories which we plan to
develop in part II of this work.

\subsection{Omega cliques}

To define the functor $\Pi $ we must first define the notion of an omega
clique and then exhibit a method for constructing omega cliques.

\begin{definition}
An \textbf{omega clique} is a bridge magma (definition 11)

$(M,R^{M})$ such that 
\[
a,b\in M_{j}\text{ and }a\parallel b\text{ }\Longrightarrow \;\exists c\text{
and }(a,b,c)\in R_{j}^{M}. 
\]
\end{definition}

Recall our convention that any two elements $a,b$ in $M_{0}$ are parallel.
Thus in an omega clique any ordered pair of 0-cells are bridged by a 1-cell
and the same is true of any ordered, parallel pair of $j$-cells for $j\geq 1$%
. \ A 1-skeletal omega clique which is also a category is a clique in the
standard sense, i.e. it is equivalent to the terminal category.

We have at hand a method for constructing omega cliques which are freely
generated by omega graphs.\ To see this let us first recall some
definitions. The category \textit{Pen\_Mor }(Section 5) has its objects
categorical Penon morphisms (definition 13). \ Its morphisms are the obvious
commutative squares. Now fix a strict omega category $A$ and let \textit{%
Pen\_Mor}$(A)$ denote the subcategory in which an object is a categorical
Penon morphism with codomain $A$ and in which a morphism has as its codomain
leg the identity functor of $A$. Let \textit{Omega\_Graph}$^{\rightarrow }$
denote the arrow category of \textit{Omega\_Graph. }Finally, let \textit{%
Omega\_Graph}$^{\rightarrow UA}$ denote the subcategory in which an object
is an arrow to the omega graph $U_{G}^{C}A$ and in which a morphism is a
commutative square in which the codomain leg is the identity on $U_{G}^{C}A$%
. In Section 7.2 we constructed an adjunction 
\[
L_{PMA}^{GMA}\dashv \;U_{GMA}^{PMA} 
\]
between the categories \textit{Omega\_Graph}$^{\rightarrow UA}$ and \textit{%
Pen\_Mor}$(A)$ in which 
\[
U_{GMA}^{PMA}:\text{\textit{Penon\_Mor}}(A)\rightarrow \text{\textit{%
Omega\_Graph}}^{\rightarrow UA} 
\]

\noindent is the obvious forgetful functor.

Now let $\Omega $ denote the terminal strict omega category. Observe that $%
U_{G}^{C}\Omega $ is the terminal omega graph. We specialize the preceding
discussion to the case $A=\Omega $.

\begin{definition}
Let $K$ be an omega graph. The \textbf{omega clique freely generated by }$K$
is the bridge magma $Cl(K)$ that is the domain of the Penon morphism $%
L_{PM\Omega }^{GM\Omega }(K\rightarrow $ $U_{G}^{C}\Omega )$.
\end{definition}

\subsection{\protect\bigskip Internal equivalences in omega magmas}

We continue laying the groundwork for the construction of the functor $\Pi $%
. \ Our next goal is to identify within any omega magma $M$ a graded subset
of cells $Eq(M)$ we shall call the \textit{internal equivalences} of $M$.

Let $K(2)$ denote the discrete omega graph with exactly two cells, denoted $%
1,2$, in dimension 0 and in which all other cells are the higher dimensional
identities associated with the cells $1,2$. \ Let $Cl(2)$ denote the free
omega clique which is the domain of the Penon morphism $L_{PM\Omega
}^{GM\Omega }(K_{2}\rightarrow $ $U_{G}^{C}\Omega )$. The omega clique $%
Cl(2) $ deserves to be called the free omega clique on two objects. For $%
i,j=1,2$ define unique 1-cells $c(i,j)$ of $Cl(2)$ by the condition 
\[
(i,j,c(i,j))\in R_{1}^{Cl(2)}\text{.} 
\]

Let $M$ be an omega magma and let $a,b$ be $(i-1)$-cells of $M$. Define the
omega magma $\mathbf{M(a,b)}$ by the following conditions:$M(a,b)_{j}\equiv
\left\{ c\in M_{i+j}\mid dom_{i-1}^{i+j}c=a\text{ and }cod_{i-1}^{i+j}c=b%
\right\} $

The functions $dom,cod$ and $id$ of $M(a,b)$ are those of $M$ suitably
restricted and reindexed

The partial operations of $M(a,b)$ are those of $M$ suitably restricted and
reindexed.

We note that any omega magma homomorphism $H:M\rightarrow N$ induces
homomorphisms 
\[
H(a,b):M(a,b)\rightarrow N(Ha,Hb) 
\]

\noindent for any pair of $(i-1)$-cells $a,b$ of $M$. \ 

\begin{definition}
Let $M$ be an omega magma. \ A $1$-cell $f$ is an \textbf{elementary
internal equivalence} if it is an identity cell or if there exists an omega
magma homorphism 
\[
H:Cl(2)\rightarrow M 
\]
such that $H(c(1,2))=f$. \ For $i\geq 2$ an $i$-cell $f$ of $M$ is an 
\textbf{elementary internal equivalence} if it is an identity cell or if
there exists an omega magma homorphism 
\[
H:Cl(2)\rightarrow M(dom_{i-2}^{i}f,\;cod_{i-2}^{i}f) 
\]
such that $H(c(1,2))=f$. Denote the graded subset of $M$ consisting of the
elementary internal equivalences by $\mathbf{el(M)}$. Note that $el(M)_{0}$
is the empty set.
\end{definition}

We note one obvious fact which follow immediately from the definition. If $M$
is an ordinary category thought of as a 1-skeletal omega magma then $%
el(M)_{1}$ ( the 1-cells of $el(M)$) consists of the isomorphisms of $M$.

Next we propose to define the concept of an internal equivalence by an
induction which starts with the elementary internal equivalences. We shall
need some additional notation. If $M$ is an omega magma and $S\subset M$ a
graded subset define $\overline{S}\subset M$ to be the smallest graded
subset of $M$ containing all composites of the elements of $S$.

\begin{definition}
Let $M$ be an omega magma and $S\subset M$ a graded subset. Define $%
T(S)\subset M$ to be the smallest graded subset which for all $k\geq 1$
contains every $k$-cell $f\in M_{k}$ satisfying the following condition:

\begin{enumerate}
\item  $\exists \;g\in M_{k}$ and $dom_{k-1}^{k}g=cod_{k-1}^{k}f$ and $%
cod_{k-1}^{k}g=dom_{k-1}^{k}f$

\item  $\exists \;h\in S$ and $dom_{k}^{k+1}h=f\bigcirc _{k-1}^{k}g$ and $%
cod_{k}^{k+1}h=id_{k}^{k-1}dom_{k-1}^{k}f$

\item  $\;h^{\prime }\in S$ and $dom_{k}^{k+1}h=g\bigcirc _{k-1}^{k}f$ and $%
cod_{k}^{k+1}h^{\prime }=id_{k}^{k-1}cod_{k-1}^{k}f$
\end{enumerate}
\end{definition}

We note that if $f,g$ are as above then $g$ is also in $T(S)$ by the
symmetry of the definition.

\begin{definition}
Define $S_{1}(M)\equiv \overline{el(M)}$ and for $i\geq 2$ define $%
S_{i}(M)\equiv \overline{S_{i-1}(M)\cup T(S_{i-1}(M))}$. \ A $k$-cell $f$ of 
$M$ is an \textbf{internal equivalence} if $f\in S_{i}(M)$ for some $i\geq 1$%
. Denote the graded subset of internal equivalences by $\mathbf{Eq(M)}$.
\end{definition}

We record the following facts.

\begin{proposition}
Let $F:M\rightarrow N$ be a morphism of omega magmas. If $f\in S_{i}(M)$
then $Ff\in S_{i}(N)$.
\end{proposition}

\noindent \textbf{Proof:}

Since $F$ is an omega magma morphism we have $F(el(M))\subseteq el(N)$. The
desired result follows by induction after observing that all the conditions
of definition 38 are preserved by omega magma homomorphisms.$\blacksquare $

\begin{proposition}
Let $a\in Eq(M)_{j}$ and define $f_{1}\equiv dom_{j-1}^{j}\alpha $ and $%
f_{2}\equiv cod_{j-1}^{j}\alpha $. Then 
\[
f_{1}\in Eq(M)_{j-1}\Longleftrightarrow f_{2}\in Eq(M)_{j-1} 
\]
\end{proposition}

\noindent \textbf{Proof:}

Assume $f_{1}\in Eq(M)_{j-1}$. Note that since $a\in Eq(M)_{j}$ there is a $%
j $-cell $\ \alpha ^{\prime }\in $ $Eq(M)_{j}$ with $dom_{j-1}^{j}\alpha
^{\prime }=f_{2}$ and $cod_{j-1}^{j}\alpha ^{\prime }=f_{1}$. \ Since $%
f_{1}\in Eq(M)_{j-1}$there exists $g\in M_{j-1}$ and $\beta ,\beta ^{\prime
}\in Eq(M)_{j}$ such that

\begin{enumerate}
\item  $dom_{j-1}^{j}\beta =f_{1}\bigcirc _{j-2}^{j-1}g$

\item  $cod_{j-1}^{j}\beta =id_{j-1}^{j-2}dom_{j-2}^{j-1}f_{1}$

\item  $dom_{j-1}^{j}\beta ^{\prime }=id_{j-1}^{j-2}dom_{j-2}^{j-1}f_{1}$

\item  $cod_{j-1}^{j}\beta ^{\prime }=f_{1}\bigcirc _{j-2}^{j-1}g$
\end{enumerate}

\noindent Then the $j$-cells $(\alpha \bigcirc
_{j-1}^{j}id_{j}^{j-1}g)\bigcirc _{j-1}^{j}\beta $

\noindent and $(\alpha ^{\prime }\bigcirc _{j-1}^{j}id_{j}^{j-1}g)\bigcirc
_{j-1}^{j}\beta ^{\prime }$are both in $Eq(M)_{j}$ \noindent and this shows
that $f_{2}\in $ $Eq(M)_{j-1}$.$\blacksquare $

\begin{proposition}
Let $\kappa :X\rightarrow Y$ be a categorical Penon morphism (definition
13). Suppose $(a.b.c)\in R_{i}^{X}$. \ Then $c\in el(X)_{i+1}$.
\end{proposition}

\noindent \textbf{Proof:}

We prove the result for $i\geq 1$; the case $i=0$ is the same but
notationally simpler.

Since $(a,b,c)\in R_{i}^{X}$ and $\kappa $ is a categorical Penon morphism
we know that $a\parallel b$, that $\kappa a=\kappa b$ and that 
\[
\kappa c=id_{i+1}^{i}a=id_{i+1}^{i}b 
\]
\noindent We simplify notation by setting $u=dom_{i-1}^{i}a$ and $%
v=cod_{i-1}^{i}a$. \noindent Note that $\kappa u=\kappa v$. \ Let $C(u,v)$
denote the omega submagma of $X(u,v)$ consisting whose $0$-cells are $a$ and 
$b$ and whose $j$-cells for $j\geq 1$ are the $j$-cells of $X(u,v)$ which $%
\kappa (u,v)$ maps to identities. Observe that the cell $c$ of $X$ is a $1$%
-cell of $C(u,v)$. \ 

We now show that $c\in Eq(X)_{i+1}$ by constructing an omega magma
homomorphism 
\[
H:Cl(2)\rightarrow C(u,v) 
\]

\noindent such that $H(c(1,2))=\dot{c}$. By definition the restriction of $%
\kappa (u,v)$ to $C(u,v)$ factors through the map from the terminal strict
omega category $\Omega $ to $Y(\kappa u,\kappa v)$ which has as its image in
dimension $0$ the cell $\kappa a$. Since $\kappa (u,v)$ is itself a
categorical Penon morphism the bridge magma $C(u,v)$ is an omega clique (its
bridge relation is the restriction of $R^{X}$ to $\ C(u,v)$ suitably
reindexed). There is an obvious map of omega graphs sending the objects $1,2$
of $K_{2}$ to $a,b\in U_{G}^{M}C(u,v)$. By adjointness we obtain a map $H$
which by definition is a bridge morphism of omega magmas and thus has the
desired property.$\blacksquare $

\subsection{Tame omega pseudo-functors}

We next identify a special subset of $Eq(M)$ which we shall term the
internal contractions of $M$.

\begin{definition}
We say that $f\in Eq(M)_{j}$ \ is an \textbf{internal contraction }and write 
$f\in Contr(M)_{j}$ if either $dom_{j-1}^{j}f$ and/or $cod_{j-1}^{j}f$ is an
identity cell. We denote the graded subset of contractions by $\mathbf{%
Contr(M)}.$
\end{definition}

\begin{definition}
Let $M$,$N$ be omega magmas and let $F:U_{S}^{M}M\rightarrow U_{S}^{M}N$ be
a morphism of the underlying globular sets. We say that $F$ is \textbf{tame}
if $F(Contr(M)\cup el(M))\subseteq Eq(N)$.
\end{definition}

Thus a morphism between the underlying globular sets of two omega magmas is
tame if it sends every internal contraction and every elementary internal
equivalence to an internal equivalence.

\begin{definition}
Let $\mathbf{F:X\rightarrow Y}$ be an omega pseudo-functor. We say that $F$
is a\textbf{\ tame omega pseudo-functor }if $F_{1}(Eq(X_{1}))\subseteq
Eq(Y_{1})$
\end{definition}

\begin{theorem}
Let $\mathbf{F:X\rightarrow Y}$ be an omega pseudo-functor. Then $\mathbf{F}$
is tame if and only if $F_{1}:U_{S}^{M}X_{1}\rightarrow U_{S}^{M}Y_{1}$ is
tame.
\end{theorem}

\noindent \textbf{Proof:}

Recall definition 40 and in particular that $S_{1}(X_{1})\equiv el(X_{1})$.
\ By hypothesis $F_{1}(S_{1}(X_{1}))\subseteq Eq(Y_{2})$. \ Now suppose $%
f\in S_{n}(X_{1})_{j}$. Then by definition there is a $g\in S_{n}(X_{1})_{j}$
and a pair of $(j+1)$-cells

\noindent\ $h,h^{\prime }\in $ $S_{n-1}(X_{1})_{j+1}$ such that

\begin{enumerate}
\item  $dom_{j}^{j+1}h=f\bigcirc _{j-1}^{j}g$ and $%
cod_{j}^{j+1}h=id_{j}^{j-1}dom_{j-1}^{j}f$

\item  $dom_{j}^{j+1}h^{\prime }=g\bigcirc _{j-1}^{j}f$ and $%
cod_{j}^{j+1}h^{\prime }=id_{j}^{j-1}cod_{j-1}^{j}f.$
\end{enumerate}

\noindent Moreover, by hypothesis we know that since $h,h^{\prime }$ are
both contractions $F_{1}(h),F_{1}(h^{\prime })\in S_{k}(Y_{1})$ for some $k$.

Next consider the cells $\rho ^{X}f$ and $\rho ^{X}g$ in $\ X_{2}$.

\noindent Denote $F_{2}\rho ^{X}f$ $\bigcirc _{j-1}^{j}F_{2}\rho ^{X}g$ by $%
u $ and $F_{2}(\rho ^{X}f$ $\bigcirc _{j-1}^{j}\rho ^{X}g)$ by $v$. Since $%
\mathbf{F}$ is an omega pseudo-functor we know that there is a $(j+1)$-cell $%
w$ such that $(u,v,w)\in R_{j}^{Y}$. By proposition 43 we know that $w\in
el(Y_{2})_{j+1}$. From proposition 41 we conclude that $\lambda ^{Y}w\in
el(Y_{1})_{j+1}$. Reversing the roles of $f,g$ in the arguments of $F_{2}$
yields a cell $\lambda ^{Y}w^{\prime }\in el(Y_{1})_{j+1}$. \ 

A very similar argument applied to $\rho ^{X}id_{j}^{j-1}dom_{j-1}^{j}f$
yields a $(j+1)$-cells $z\in el(Y_{2})_{j+1}$ whose domain is $F_{2}(\rho
^{X}(id_{j}^{j-1}dom_{j-1}^{j}f))$ and whose codomain is $%
id_{j}^{j-1}(F_{2}(dom_{j-1}^{j}\rho ^{X}f))$ while replacing $f$ by $g$
yields a similar cell $z^{\prime }$.

We complete the proof by considering the $(j+1)$-cells

\noindent $(\lambda ^{Y}w\bigcirc _{j}^{j+1}F_{1}h)\bigcirc
_{j}^{j+1}\lambda ^{Y}z$ and $(\lambda ^{Y}z^{\prime }\bigcirc
_{j}^{j+1}F_{1}h^{\prime })\bigcirc _{j}^{j+1}\lambda ^{Y}w^{\prime }$.Each
of these cells is an element of $S_{k}(Y_{1})$ and hence their composites
are elements of $S_{k+1}(Y_{1})$. \ The first of these two cells has its
domain the $j$-cell $(F_{1}f\bigcirc _{j-1}^{j}F_{1}g)$ and its codomain the 
$j$-cell $(id_{j}^{j-1}dom_{j-1}^{j}f)$. The second has as its domain $%
(F_{1}g\bigcirc _{j-1}^{j}F_{1}f)$ and as its codomain $%
(id_{j}^{j-1}cod_{j-1}^{j}f\dot{)}$. This shows that $F_{1}f\in
S_{k+2}(Y_{1})$. $\blacksquare $

\begin{corollary}
Let $\mathbf{F:X\rightarrow Y}$ be a tame omega pseudo-functor. Then $%
\mathbf{F(}a,b)$ is also a tame omega pseudo-functor for any choice of $%
(i-1) $-cells $a,b.$
\end{corollary}

The preceding theorem naturally leads one to ask which classes of omega
pseudo-functors can a priori be identified as tame. By proposition 41 any
omega functor $\mathbf{F}$ is a tame omega pseudo-functor because, by
definition, $F_{1}$ is a morphism of omega magmas.

\begin{theorem}
Let $\mathbf{F:X\rightarrow Y}$ be an omega pseudo-functor. If both $\mathbf{%
X}$ and $\mathbf{Y}$ are n-skeletal then $\mathbf{F}$ is tame.
\end{theorem}

\noindent \textbf{Proof:}

We show by a downward induction on the dimension of cells in $Eq(X_{1})$
that $F_{1}(Eq(X_{1}))\subseteq Eq(Y_{1})$.

Let $f$ be an $n$-cell of $X_{1}$ and suppose $f\in Eq(X_{1})$. $\ $Note
that $f$ is a 1-cell in $X(dom_{n-1}^{n}f,cod_{n-2}^{n}f)_{1}$. Now $\mathbf{%
X}(dom_{n-1}^{n}f,cod_{n-2}^{n}f)$ is an 1-skeletal weak omega category and
hence $X(dom_{n-1}^{n}f,cod_{n-2}^{n}f)_{1}$ is a category by the results of
Section 3. If $A$ is a category (considered as a 1-skeletal omega magma) it
follows immediately from definitions 39 and 40 that $Eq(A)_{1}$ is the set
of isomorphisms of $A$. \ Now $F(dom_{n-1}^{n}f,cod_{n-2}^{n}f)_{1}$ (the
first coordinate of $\mathbf{F}(dom_{n-1}^{n}f,cod_{n-2}^{n}f)$) is a
functor and functors preserve isomorphisms. Thus $F_{1}f\in Eq(Y_{1})_{n}$
and thus we conclude that $F_{1}(Eq(X_{1})_{n})$ $\subseteq Eq(Y_{1})_{n}$.

Now suppose that $F_{1}(Eq(X_{1})_{j})\subseteq Eq(Y_{1})_{j}$ for all $%
j\geq k+1$. \ Let $f\in Eq(X_{1})_{k}$ and let $g$ be the $k$-cell and $%
h,h^{\prime }$ be the $(k+1)$-cells which must witness this property as in
definition 39. As in the proof of the preceding theorem we can produce a
composite $(k+1)$-cell in $Y_{1}$ with domain $(F_{1}g\bigcirc
_{j-1}^{j}F_{1}f)$ and codomain $(id_{j}^{j-1}cod_{j-1}^{j}f\dot{)}$. Two of
the three $(k+1)$-cells making up this composite are in $el(Y$ $_{1})$ and
the third, $F_{1}h$, is in $Eq(Y_{1})_{k+1}$ by induction. Thus the
composite is in $Eq(Y_{1})_{k+1}$. This coupled with the same argument with $%
h^{\prime }$ in place of $h$ shows that $F_{1}f\in Eq(Y_{1})_{k+1}$ and thus
completes the inductive step.$\blacksquare $

Next we delineate another class of tame, omega pseudo-functors.

\begin{theorem}
Let $\mathbf{X}$ be a weak omega category, $a,b,c$ three $(i-1)$-cells of $%
X_{1}$ and $h$ an $i$-cell of $X_{1}$ such that $dom_{i-1}^{i}h=b$ and $%
cod_{i-1}^{i}h=c$. Then $\mathbf{\Theta (h):X(}a,b)\rightarrow $ $\mathbf{X(}%
a,c)$ (Section 8.2) is a tame, omega pseudo-functor.
\end{theorem}

\noindent \textbf{Proof:}

$Eq(X_{1})$ is closed under composition and contains all the identities of $%
X_{1}$. The result then follows from the definition of $\Theta (h)_{1}$
which is just composition with the various identities of $h$. $\blacksquare $

\subsection{The functor $\Pi $}

Let $\mathbf{X}$ be a weak omega category. \ We define an equivalence
relation $\approx $ on the $i$-cells of $X_{1}$ by setting 
\[
a\approx b\text{ }\Longleftrightarrow \;\exists f\in Eq(X_{1})_{i+1}\text{
and }dom_{i}^{i+1}f=a\text{ and }cod_{i}^{i+1}f=b 
\]
Define $\Pi \mathbf{X}$ to be the set of $\approx $ equivalence classes of $%
(X_{1})_{0}$, the $0$-cells of $X_{1}$.

Let \textit{Tame\_Omega\_Cat }denote the category in which an object is a
small, weak omega category and in which a morphism is a tame pseudo-functor.
It contains \textit{Omega\_Cat }as a subcategory and, by theorem 49, it
contain as a full subcategory the full subcategory of \textit{PS\_Omega\_Cat}
(in which objects are weak omega categories and morphisms arbitrary omega
pseudo-functors) containing all n-skeletal, weak omega categories for all $n$%
.

Theorem 47 then yields:

\begin{theorem}
$\Pi :$\textit{Tame\_}Omega\_Cat$\rightarrow $Set is a functor.
\end{theorem}

\noindent

\subsection{Weak equivalence and omega equivalence}

Let $\mathbf{F:X\rightarrow Y}$ be a tame omega pseudo-functor and for any
pair of $(i-1)$-cells $a,b$ of $X_{1}$ let 
\[
\mathbf{F}(a,b):\mathbf{X}(a,b)\rightarrow \mathbf{Y}(F_{1}a,F_{1}b) 
\]

\noindent be the tame omega pseudo-functor associated with $\mathbf{F},a,b$.

\begin{definition}
$\mathbf{F}$ is a \textbf{weak equivalence} if $\Pi \mathbf{F}$ and $\Pi 
\mathbf{F}(a,b)$ are isomorphisms for all $(i-1)$-cells $a,b$ and all $i\geq
1$.
\end{definition}

\begin{definition}
$\mathbf{F}$ is an \textbf{omega equivalence} if there is a tame omega
pseudo-functor $\mathbf{G:Y\rightarrow X}$ such that 
\[
\Pi (\mathbf{G\circ F),\;}\Pi (\mathbf{F\circ G)},\noindent \;\Pi \mathbf{%
(G\circ F)}(a,b),\;\Pi \mathbf{(F\circ G)}(a,b) 
\]
are isomorphisms for all $(i-1)$-cells $a,b$ and all $i\geq 1$.
\end{definition}

\begin{theorem}
Let $\mathbf{F:X\rightarrow Y}$ be an omega pseudo-functor. Suppose that $%
\mathbf{F}$ is a weak equivalence and that $\mathbf{X,Y}$ are 1-skeletal.
Then $F_{1}$ is an equivalence of ordinary categories.
\end{theorem}

\noindent \textbf{Proof:}

Since $\mathbf{X,Y}$ are 1-skeletal we know that $\mathbf{F}$ is tame and
from Section 3 we know that $F_{1}$ is a functor between ordinary
categories. We have already observed that in this situation $Eq(X_{1})_{1}$
consists of the isomorphisms of the category $X_{1}$. Thus the functor $%
F_{1} $ is essentially surjective since it induces a bijection on
isomorphism classes. The fact that $X_{1}$ is 1-skeletal means that $%
Eq(X_{1})_{j}$ consists only of identities for $j\geq 2$ so $F_{1}$ must be
fully faithful. $\blacksquare $

\begin{theorem}
Let $\mathbf{F:X\rightarrow Y}$ be an omega pseudo-functor. Suppose that $%
\mathbf{F}$ is a weak equivalence and that $\mathbf{X,Y}$ are 2-skeletal.
Then $F_{1}$ is a biequivalence between bicategories.
\end{theorem}

\noindent \textbf{Proof:}

Since $\mathbf{X,Y}$ are 2-skeletal we know that $\mathbf{F}$ is tame and
from the results of Sections 4 and 8 we know that $F_{1}$ is a classical
pseudo-functor between bicategories.\ For any pair of 0-cells $a,b$ of $%
X_{1} $ the preceding theorem says that the functor 
\[
\mathbf{F}(a,b):\mathbf{X}(a,b)\rightarrow \mathbf{Y}(F_{1}a,F_{1}b) 
\]

\noindent is an equivalence of ordinary categories. The fact that $F_{1}$ is
essentially surjective follows from the observation that the elements of $%
Eq(X_{1})$ are precisely the $1$-cell equivalences in the ordinary,
bicategorical sense. $\blacksquare $

We next exhibit the most common examples of omega equivalences between weak
omega categories.

\begin{proposition}
Let $\mathbf{X}$ be a weak omega category and let $a,b$ be a pair of $(i-1)$%
-cells of $X_{1}$. Then the omega pseudo-functor 
\[
\mathbf{\Theta (id}_{i}^{i-1}\mathbf{b}):\mathbf{X}(a,b)\rightarrow \mathbf{X%
}(a,b)\noindent \noindent 
\]
is an omega equivalence. In fact $\Pi (\mathbf{\Theta (id}_{i}^{i-1}\mathbf{b%
}))$ is the identity on $\Pi \mathbf{X}(a,b)$.
\end{proposition}

\noindent \textbf{Proof:}

By theorem 50 $\mathbf{\Theta (id}_{i}^{i-1}\mathbf{b})$ is a tame omega
pseudo-functor. That it is an omega equivalence follows immediately from
propositions 43 and 41 and the fact that therefore \ 
\[
f\bigcirc _{i-1}^{i}id_{i}^{i-1}b\approx f 
\]
\ via an internal equivalence (which in fact is an elementary internal
equivalence). This observation also shows that $\Pi (\mathbf{\Theta (id}%
_{i}^{i-1}\mathbf{b}))$ is the identity on $\Pi \mathbf{X}(a,b)\blacksquare $

\begin{proposition}
Let $\mathbf{X}$ be a weak omega category and let $a,b,c$ be three $(i-1)$%
-cells of $X_{1}$. Let $h_{1},h_{2}\,$ be $i$-cells with domain $b$ and
codomain $c$. Suppose $h_{1}\approx h_{2}$. Then $\mathbf{\Theta (h}_{1}%
\mathbf{)}$ is an omega equivalence if and only if $\mathbf{\Theta (h}_{1}%
\mathbf{)}$ is an omega equivalence.
\end{proposition}

\noindent \textbf{Proof:}

This again is a straightforward consequence of the fact that $Eq(X_{1})$
contains all the identity cells of $X_{1}$ and is closed under composition. $%
\blacksquare $

\begin{theorem}
Let $\mathbf{X}$ be a weak omega category and let $a,b,c$ be three $(i-1)$%
-cells of $X_{1}$. Let $h$ be an $i$-cell with domain $b$ and codomain $c$.
If $h\in Eq(X_{1})$ then $\mathbf{\Theta (h)}$ is an omega equivalence.
\end{theorem}

\noindent \textbf{Proof:}

If $h\in Eq(X_{1})$ then by definition there exists $h^{\prime }\in
Eq(X_{1}) $ such that 
\[
h\bigcirc _{i-1}^{i}h^{\prime }\approx id_{i}^{i-1}b 
\]
\[
h^{\prime }\bigcirc _{i-1}^{i}h\approx id_{i}^{i-1}c 
\]

The desired result now follows from the preceding two propositions and the
fact that $\Pi $ is a functor.$\blacksquare $

\end{document}